\documentclass[12pt]{article}

\usepackage{amsfonts,amsmath,amssymb,bbm}
\usepackage{graphicx}
\usepackage{eufrak}

\usepackage{setspace} 
\def\proof{\noindent{\bf Proof:}\hskip10pt}
\def\proofof{\noindent{\bf Proof of }}        
\def\QED{\hfill $\Box$}

\font\tenmath=msbm10 scaled 1200
\font\sevenmath=msbm7 scaled 1200
\font\Fivemath=msbm5 scaled 1200
\newfam\mathfam \textfont\mathfam=\tenmath
\scriptfont\mathfam=\sevenmath \scriptscriptfont\mathfam=\Fivemath

\def \\ { \cr }

\def \1{\mathbbmss{1}}
\def\N{\mathbb{N}}
\def\E{\mathbb{E}}
\def\P{\mathbb{P}}

\def \e{{\rm e}}
\def \ebf{{\rm \bf e}}

\def \ve{{\ensuremath{\varepsilon}}}
\def \dt{{\ensuremath{\textup{d}}}}
\def \pB{{\ensuremath{\mathcal{B}}}}
\def \pF{{\ensuremath{\mathcal{F}}}}
\def \pG{{\ensuremath{\mathcal{G}}}}

\def \pU{{\ensuremath{\mathcal{U}}}}

\def \pZ{{\ensuremath{\mathcal{Z}}}}
\def \pC{{\ensuremath{\mathcal{C}}}}
\def \pCf{{\ensuremath{\mathfrak{C}}}}

\newcommand{\eql}{\ensuremath{\overset{\mbox{\scriptsize\textup{d}}}{=}}}
\newtheorem{theorem}{Theorem}

\newtheorem{proposition}{Proposition}
\newtheorem{lemma}{Lemma}
\newtheorem{corollary}{Corollary}
\usepackage[hyperfootnotes=false]{hyperref}

\begin{document}

\title{Percolation on random recursive trees}\author{Erich
  Baur\footnote{UMPA ENS de Lyon, 46, allée d'Italie, 69364 Lyon Cedex 07, France. Email:
  erich.baur@math.uzh.ch}\\{ENS Lyon}}

\maketitle 
\thispagestyle{empty}

\begin{abstract} We study Bernoulli bond percolation on a random
  recursive tree of size $n$ with percolation parameter $p(n)$ converging
  to $1$ as $n$ tends to infinity. The sizes of the percolation clusters
  are naturally stored in a tree. We prove convergence in distribution of
  this tree to the genealogical tree of a continuous-state branching
  process in discrete time. As a corollary we obtain the asymptotic sizes
  of the largest and next largest percolation clusters, extending thereby
  a recent work of Bertoin~\cite{Be1} which deals with cluster sizes in the
  supercritical regime. In a second part, we show that the
  same limit tree appears in the study of the tree components which emerge
  from a continuous-time destruction of a random recursive tree. We comment
  on the connection to our first result on Bernoulli bond percolation.
\end{abstract} 

\footnote{{\it Acknowledgment of support.} This
  research was supported by the Swiss National Science Foundation grant
  P2ZHP2\_151640.}

\section{Introduction}
\label{Sintro}
Let $V$ be a finite and totally ordered set of vertices. An increasing tree
on $V$ is a tree rooted at the smallest element of $V$ such that the
sequence of vertices along the branch from the root to any vertex
increases.  A {\it random recursive tree} (RRT for short) of size $n+1$ is
a tree picked uniformly at random amongst all increasing trees on
$\{0,1,\dots,n\}$. Henceforth we write $T_n$ for such a RRT. Note
that the root vertex of $T_n$ is given by $0$.

We consider Bernoulli bond percolation on $T_n$ with parameter $p(n)\in
(0,1)$. This means we first pick $T_n$ and then remove each edge with
probability $1-p(n)$, independently of the other edges. We obtain a
partition of vertices into clusters, i.e. connected components, and we are
concerned with the asymptotic sizes of these clusters. Let us call
percolation on a RRT $T_n$ in the regime
$p(n)\rightarrow 1$
\begin{itemize}
\item {\it weakly supercritical}, if $\frac{1}{\ln n}\ll 1-p(n)\ll 1$,
\item {\it supercritical}, if $1-p(n)\sim \frac{t}{\ln n}$ for some $t>0$
  fixed,
\item {\it strongly supercritical}, if $0<1-p(n)\ll \frac{1}{\ln n}$.
\end{itemize}
The terminology is explained by our results: We will see that the root
cluster has size $\sim n^{p(n)}$, while the next largest
clusters have a size of order $(1-p(n))n^{p(n)}$.

We encode the sizes of all percolation clusters by a tree structure, which
we call the {\it tree of cluster sizes}.  A percolation cluster of $T_n$ is
called a cluster of generation $k$, if it is disconnected from the root
cluster by exactly $k$ deleted edges. In the tree of cluster sizes,
vertices of level $k$ are labeled by the sizes of the clusters of
generation $k$. Consequently, the root vertex represents the size of the
root cluster of $T_n$.  Then, if a vertex represents the size of a cluster
$\tau$ of generation $k$, its children are given by the sizes of those
clusters of generation $k+1$ which are separated from $\tau$ by one deleted
edge (see Figure $1$).

We normalize cluster sizes of generation $k$ by a factor
$(1-p(n))^{-k}n^{-p(n)}$. After a local re-ranking of vertices, we show
that the tree of cluster sizes converges in distribution to the
genealogical tree of a continuous-state branching process (CSBP) in
discrete time, with reproduction measure $\nu(da) = a^{-2}\dt a$ on
$(0,\infty)$ and started from a single particle of size $1$ (Theorem
\ref{T1}). Moreover, we obtain precise limits for the largest non-root
clusters (Corollary \ref{C1}).

Asymptotic cluster sizes have been studied for numerous other random graph
models. At first place, these include the Erd\"os-R\'enyi graph
model (see Alon and Spencer \cite[Chapter 11]{AS} for an overview with
further references). Concerning trees, uniform Cayley trees of size $n$
have been studied in the works of Pitman \cite{Pi}, \cite{Pi2} and Pavlov
\cite{Pa} in the regime $1-p(n)\sim t/\sqrt{n}$, $t>0$ fixed, where the
number of giant components is unbounded. For general large
trees, Bertoin gives in \cite{Be2} a criterion for the root cluster of a
Bernoulli bond percolation to be the (unique) giant cluster.

Unlike random Cayley trees of size $n$, whose heights are typically of
order $\sqrt{n}$, random recursive trees have heights of logarithmic order
(see e.g. the book of Drmota \cite{Dr}). Bertoin proved in \cite{Be1} that
in the supercritical regime when $1-p(n)\sim t/\ln n$, the size of the root
cluster of a RRT on $n+1$ vertices, normalized by a factor $1/n$, converges
to $\e^{-t}$ in probability, while the sizes of the next largest clusters,
normalized by a factor $\ln n/ n$, converge to the atoms of some Poisson
random measure. This result was extended by Bertoin and Bravo \cite{BeBr}
to large scale-free random trees, which grow according to a preferential
attachment algorithm and form another family of trees with
logarithmic height.

Here we follow the route of \cite{BeBr} and analyze a system of branching
processes with rare neutral mutations. In this way we gain control over the
sizes of the root cluster and of the largest clusters of the first
generation, for all regimes $p(n)\rightarrow 1$. An iteration of the
arguments then allows us to prove convergence of higher generation cluster
sizes.

The methods of \cite{Be1} are based on a coupling of Iksanov and M\"ohle
\cite{IM} between the process of isolating the root in a RRT and a certain
random walk. They seem less suitable for the weakly supercritical regime,
where one has to look beyond the passage time up to which the coupling is
valid. This was already mentioned in the introduction of \cite{Be1}, where
also the question is raised how the sizes of the largest clusters behave
when $1-p(n)\gg 1/\ln n$. 

In the second part of this paper, we however
readopt the methods of \cite{Be1}. We consider a destruction process on
$T_n$, where edges are equipped with i.i.d. exponential clocks and deleted
at the time given by the corresponding variable. Starting with the full
tree $T_n$, each removal of an edge $e$ gives birth to a new tree component
rooted at the outer endpoint of $e$. The order in which the tree components
are cut suggests an encoding of their sizes and birth times by a
tree-indexed process, which we call the {\it tree of components} (see
Figure $2$ for an example).

Keeping track of the birth times allows us to consider only those tree
components which are born in the destruction process up to a certain
finite time. Interpreting the latter as a version of a Bernoulli bond
percolation on $T_n$, tree components are naturally related to percolation
clusters. This observation was made by Bertoin in \cite{Be1} and then used
to study cluster sizes in the supercritical regime. We further develop
these ideas in the last section and make the link to our results on
percolation from the first part of this paper. The starting point is a
limit result for the tree of components (Theorem \ref{T2}), which we
believe is of interest on its own.

The destruction process can be viewed as an iterative application of the
cutting down or isolation of the root process, which has been analyzed in
detail for RRT's in Meir and Moon \cite{MM2}, Panholzer, Drmota {\it et
  al.} \cite{DIMR}, Iksanov and M\"ohle \cite{IM}, Bertoin \cite{Be1} and
others. The tree of components should be seen as a complement to the
so-called {\it cut-tree}, which was studied for random recursive trees by
Bertoin in~\cite{Be3}. We briefly recall its definition at the very end.

The rest of this paper is organized as follows. The goal of Section
\ref{Streeofclustersizes} is to prove our main result Theorem \ref{T1} on
the sizes of percolation clusters. We first introduce the tree of cluster
sizes and state the theorem. Then we establish the connection to Yule
processes and obtain the asymptotic sizes of the root and first generation
clusters. We then turn to higher generation clusters in the tree and finish
the proof of Theorem \ref{T1}. Section \ref{Streeofcomponents} is devoted
to the analysis of the destruction process of a RRT. At first we define the
tree of components and formulate our main result Theorem \ref{T2} for this
tree. The splitting property of random recursive trees transfers into a
branching property for the tree of components, which we illustrate together
with the coupling of Iksanov and M\"ohle in the next part. Then we prove
Theorem \ref{T2}. In the last part, we sketch how our analysis
of the destruction process leads to information on percolation clusters in
the supercritical regime. We compare our results there with
Theorem \ref{T1} and finish this paper by pointing at the connection to the
cut-tree. 

We finally mention that except for the very last part, Section
\ref{Streeofcomponents} on the destruction process can be
read independently of Section \ref{Streeofclustersizes}.

\section{The tree of cluster sizes}
\label{Streeofclustersizes}
\subsection{Our main results on percolation}
\label{Smainresults}
We use a tree structure to store the percolation clusters or, more
precisely, their sizes. In this direction, recall that the universal tree
is given by
\begin{equation*}
\pU = \bigcup_{k=0}^\infty \N^k,
\end{equation*}
with the convention $\N^0 = \{\emptyset\}$ and $\N=\{1,2,\dots\}$. In
particular, an element $u\in\pU$ is a finite sequence of strictly positive
integers $(u_1,\dots,u_k)$, and we refer 
to its length $|u|=k$ as the ``generation'' or level of $u$. The $j$th child
of $u$ is given by $uj=(u_1,\dots,u_k,j)$, $j\in\mathbb{N}$. The empty
sequence $\emptyset$ is the root of the tree and has length
$|\emptyset|=0$. If no confusion occurs, we simply write $u_1\dots u_k$ instead of
$(u_1,\dots,u_k)$. 

Now consider Bernoulli bond percolation on a RRT $T_n$ with parameter
$p(n)$. This induces a family of percolation clusters, and we say that a
cluster is of generation $k=0,1,\dots,n$, if it is disconnected from the
root $0$ by $k$ erased edges. This means that exactly $k$ edges have been
removed by the percolation from the path in the original tree $T_n$
connecting $0$ to the root of the cluster. In this terminology, the only
cluster of the zeroth generation is the root cluster.

We define recursively a process $\pC^{(n)} =(\pC_u^{(n)}:u\in\pU)$
indexed by the universal tree, which we call the {\it tree of cluster
  sizes}. 

First, $\pC_\emptyset^{(n)}$ is the size of the root cluster of
$T_n$. Next, we let $\pC_1^{(n)}\geq \pC_2^{(n)} \geq \dots\geq \pC_{\ell(n)}^{(n)}$
denote the decreasingly ranked sequence of the sizes of the clusters of
generation $1$, where $\ell(n)\leq n$, and in the case of ties, clusters of
the same size are ordered uniformly at random. We continue the definition
iteratively as follows. Assume that for some $u\in\pU$ with $1\leq |u|\leq
n-1$, $\pC_u^{(n)}$ has already been defined to be the size of some cluster
$c_u^{(n)}$ of generation $|u|$. We then specify the children of
$\pC_u^{(n)}$. Among all clusters of generation $|u|+1$, we consider those
which are disconnected by exactly one erased edge from
$c^{(n)}_u$. Similar to above, we rank these clusters in the decreasing
order of their sizes and let $\pC_{uj}^{(n)}$ be the size of the $j$th
largest. An example is given in Figure $1$.

The definition is completed by putting $\pC_{u}^{(n)}=0$ for all
$\pC_u^{(n)}$ which have not been specified in the above way. In
particular, $\pC_u^{(n)}=0$ for all $u$ with $|u|>n$, and if $\pC_u^{(n)} =
0$ for some $u$, then all elements of the subtree of $\pC^{(n)}$ rooted at $u$ are set to
zero. 
\begin{figure}
\begin{center}\parbox{7cm}{\includegraphics[width=5cm]{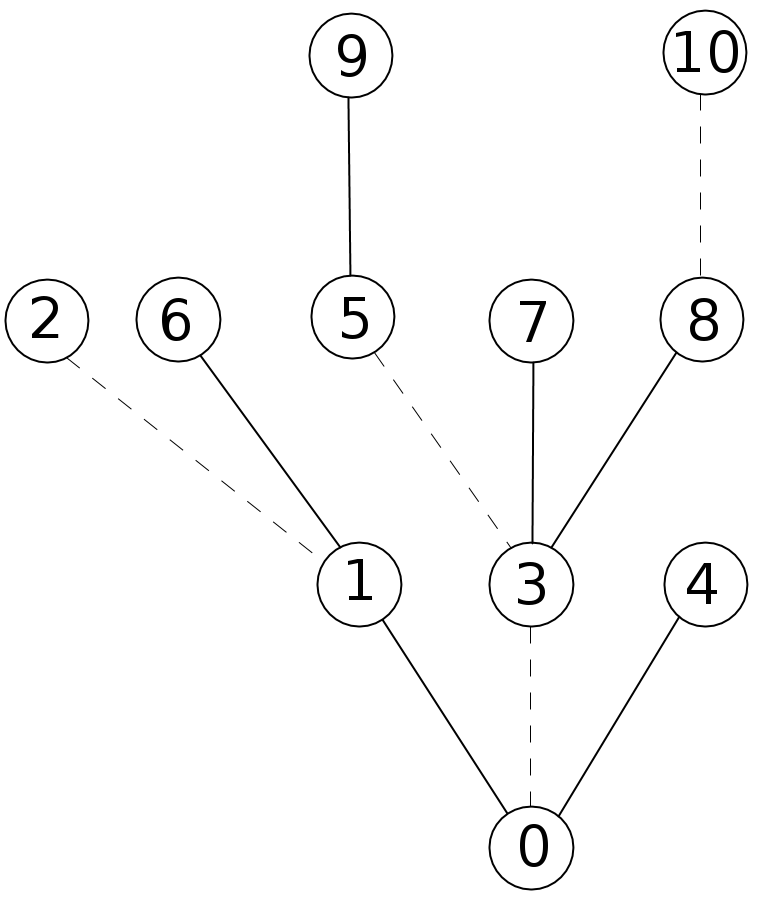}}
\parbox{7cm}{\includegraphics[width=7cm]{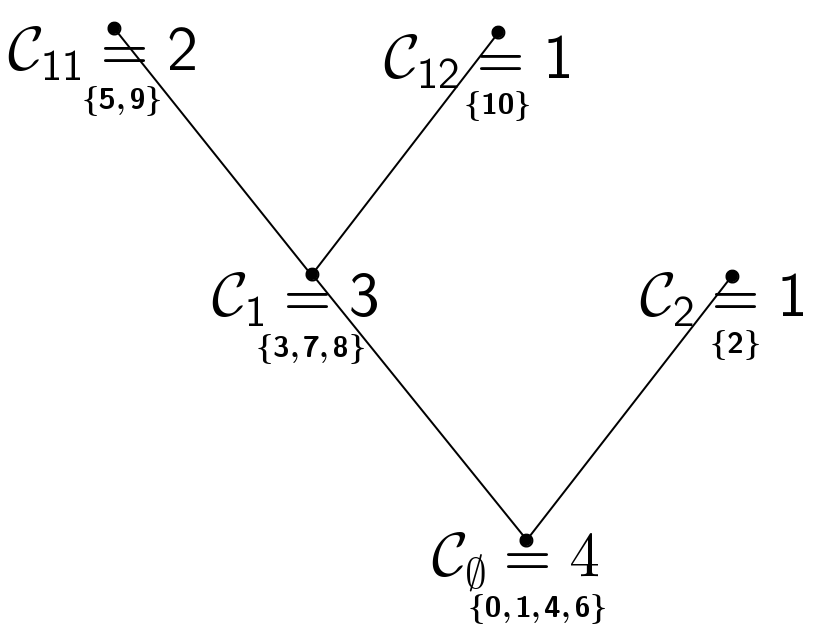}}
\end{center}
\centerline{\bf Figure 1}
\centering{\sl Left: Percolation on a recursive tree with vertices labeled
  $0,1,\dots,10$. The
  edges which were removed by the percolation are indicated by dashed lines.}\\
{\sl  Right: The corresponding percolation clusters, whose sizes are encoded by
$\pC^{(n)}$.}
\end{figure}

Our limit object is given by the genealogical tree $\pZ=(\pZ_u:u\in \pU)$
of a continuous-state branching process in discrete time with reproduction
measure $\nu(da)=a^{-2}\dt a$, started from a single particle. The
distribution of $\pZ$ is characterized by induction on the generations as
follows (cf. \cite[Definition 1]{Be}).
\begin{enumerate}
\item $\pZ_\emptyset = 1$ almost surely;
\item for every $k=0,1,2,\dots,$ conditionally on $(\pZ_v: v\in\pU,|v|\leq
  k)$, the sequences $(\pZ_{uj})_{j\in\N}$ for the vertices $u\in\pU$ at
  generation $|u|=k$ are independent, and each sequence
  $(\pZ_{uj})_{j\in\N}$ is distributed as the family of the atoms of a Poisson
  random measure on $(0,\infty)$ with intensity $\pZ_u\nu$, where the atoms
  are ranked in the decreasing order of their sizes.
 \end{enumerate}
 We turn now to the statement of Theorem \ref{T1}. Recall that we consider
 the regimes $p(n)\in(0,1)$ with $p(n)\rightarrow 1$ as
 $n\rightarrow\infty$. In the strongly supercritical regime when $1-p(n)\ll
 1/\ln n$, the root cluster has size $\sim n$, and if $(1-p(n))n^{p(n)}$
 stays bounded, the next largest clusters will be of constant size only. In
 order to exclude this case, and since we would like to consider also
 higher generation clusters, we shall implicitly assume that
 $(1-p(n))^{k}n^{p(n)}\rightarrow \infty$ for every $k\in\mathbb{N}$. If
 this last condition fails, then our convergence results do still hold
 restricted to generations $k\leq \max\{\ell\in\mathbb{N}_0 :
 (1-p(n))^{\ell}n^{p(n)}\rightarrow \infty \}$.
\begin{theorem}\label{T1} 
  As $n\rightarrow\infty$, in the sense of finite-dimensional
  distributions,
$$
\left(\frac{(1-p(n))^{-|u|}}{n^{p(n)}}\pC_u^{(n)} :u\in\pU\right)\Longrightarrow
(\pZ_u : u\in\pU).$$
\end{theorem}
While the theorem specifies the size $\pC^{(n)}_{\emptyset}$ of the
root cluster as $n$ tends to infinity, it does not immediately
answer the question how the sizes
$\mathbf{C}^{(n)}_1\geq\mathbf{C}^{(n)}_2\geq\dots$ of the next largest
clusters behave.  We will however see that for fixed $\ell$, the $\ell$
largest non-root clusters are with high probability to be found amongst the
$k$ largest clusters of the first generation, provided $n$ and $k$ are
taken sufficiently large.

\begin{corollary}\label{C1}
  As $n\rightarrow\infty$, we have $\pC^{(n)}_{\emptyset}\sim n^{p(n)}$ in
  probability. Moreover, for each fixed $\ell\in\mathbb{N}$, 
  $$
  \left(\frac{(1-p(n))^{-1}}{n^{p(n)}}\mathbf{C}^{(n)}_1,\dots,\frac{(1-p(n))^{-1}}{n^{p(n)}}\mathbf{C}^{(n)}_\ell\right)\Longrightarrow
  (\pZ_1,\dots,\pZ_\ell),$$ where in accordance with our definition of $\pZ$, $\pZ_1> \pZ_2>\dots$ are the atoms of a Poisson
  random measure on $(0,\infty)$ with intensity $\nu(da) = a^{-2}\dt a$.
\end{corollary}

\subsection{Connection to Yule processes with neutral mutations }
\label{SYule}
Here we will develop the methods that enable us to prove
finite-dimensional convergence of the tree of cluster sizes for generations
$\leq 1$. We conclude this part with the proof of Corollary \ref{C1}. In the next section, we lift the
convergence to higher levels in the tree and thereby finish the proof
Theorem \ref{T1}.

The following recursive construction of a RRT forms the basis of our
approach. We consider a standard Yule process $Z = (Z(t) : t\geq 0)$, i.e a
continuous-time pure birth process started from $Z(0) = 1$, with unit birth
rate per unit population size. Then, if the ancestor is labeled by $0$ and
the next individuals are labeled in the increasing order of their birth
times, the genealogical tree of the Yule process stopped at the instant
$$\rho_n = \inf\{t\geq 0: Z(t) = n+1\}$$ is a version of $T_n$.

With this construction, percolation on a RRT $T_n$ with parameter
$p\in(0,1)$ can be interpreted in terms of neutral mutations which are
superposed to the genealogical tree. In the description that follows we are guided
by \cite{Be} and \cite{BeBr}.

Except for the ancestor, we let each individual of the Yule process be a
clone of its parent with probability $p$ and a mutant with probability
$1-p$. Being a mutant means that the individual receives a new genetic type
which was not present before. The reproduction law is neutral in the sense
that it is not affected by the mutations. We record the genealogy of types
by the universal tree in the following way. Every vertex $u\in\pU$ stands
for a new genetic type. The empty set $\emptyset$ represents the type of
the ancestor, and for every $u=(u_1,\dots, u_k)\in\pU$ and $j\in\N$, the
$j$th child of $u$, i.e. $uj = (u_1,\dots,u_k,j)$, stands for the genetic
type which appeared at the instant when the $j$th mutant was born in the
subpopulation of type $u$.

Starting from $Z_\emptyset^{(p)}(0) = 1$ and $Z_u^{(p)}(0) = 0$ for
$u\in\pU\backslash \emptyset$, we write $Z_u^{(p)}(t)$ for the size of the
subpopulation of type $u$ at time $t\geq 0$, when neutral mutations occur
at rate $1-p$ per unit population size. Clearly, the sum over all
subpopulations $Z=(\sum_{u\in\pU}Z_u^{(p)}(t): t\geq 0)$ evolves as a standard Yule
process, and we will henceforth work with $Z$ defined in this way.

Moreover, interpreting the genealogical tree of $Z(\rho_n)$ as a RRT $T_n$
as above, the sizes of the clusters of generation $k$ are given by the
variables $Z_u^{(p)}(\rho_n)$ with $|u|=k$. Note however that in the tree
of cluster sizes, the children of each element are decreasingly ordered
according to their sizes, while in the population model, the sequence
$(Z_{uj}^{(p)}(\rho_n) : j\in\mathbb{N})$ for $u\in\pU$ is ordered
according to the birth times of the mutants stemming from type $u$,
i.e. type $ui$ was born before type $uj$ for $i<j$.

Let us denote the birth time of the subpopulation of type $u\in\pU$ by
$$
b_u^{(p)}=\inf\{t\geq 0: Z_u^{(p)}(t)>0\}.
$$
Clearly, for $p<1$, each variable $b_u^{(p)}$ is almost surely
finite. Moreover, each process $(Z_u^{(p)}(b_u^{(p)}+t) : t\geq 0)$ for
$u\in\pU$ is distributed as a continuous-time pure birth process with birth
rate $p$ per unit population size, started from a single particle. Once an
individual of a new genetic type appears, the population of that type
evolves independently, which shows that the processes
$Z_u^{(p)}(b_u^{(p)}+\cdot)$ for $u\in\pU$ are independent. The sequence of
subpopulations bearing a single mutation is moreover independent from the
sequence of its birth times:
\begin{quote}
The processes $(Z_i^{(p)}(b_i^{(p)}+t) : t\geq
0)$ for $i\in\N$ are i.i.d. and independent of the sequence of birth times
$(b_i^{(p)} : i\in \N)$.
\end{quote}
For a formal proof of this statement, see \cite[Lemma 1]{BeBr}.

Our first aim is to obtain a joint limit law for
$Z_\emptyset^{(p)}(\rho_n)$ and $Z_i^{(p)}(\rho_n)$, $i\in\N$, when
$n\rightarrow \infty$ and $p=p(n)\rightarrow 1$. In this direction, we
recall that $$W(t) = \e^{-t}Z(t),\quad t\geq 0,$$ is non-negative
square-integrable martingale with terminal value $W(\infty)$ given by a
standard exponential variable. The next lemma, which is similar to Lemma 2
in \cite{BeBr}, shows that the speed of convergence is exponential.

\begin{lemma}
\label{L2BeBr}
For every $t\geq 0$, one has 
$$\E\left(\sup_{s\geq t}|W(s)-W(\infty)|^2\right)\leq 10
\e^{-t},\quad\hbox{and}\quad\E\left(\sup_{s\geq 0}\e^{2s/3}|W(s)-W(\infty)|^2\right)\leq 
\frac{10\e^{2/3}}{(1-\e^{-1/6})^2}.
$$
\end{lemma}
Similarly, for $p\in(0,1)$ fixed, $W_\emptyset^{(p)}(t) = \e^{-pt}Z_\emptyset^{(p)}(t)$, $t\geq 0$,
is a martingale. Its terminal value $W_{\emptyset}^{(p)}(\infty)$ is
another standard exponential random variable, but for $p$ tending to $1$,
$W_{\emptyset}^{(p)}(\infty)$ converges to $W(\infty)$. More specifically,
\cite[Lemma 3]{BeBr} reads in our case as follows.
\begin{lemma}
\label{L3BeBr}
$$\lim_{p\rightarrow 1,\,t\rightarrow\infty}\E\left(\sup_{s\geq t}\left|W_{\emptyset}^{(p)}(s)-W(\infty)\right|^2\right)=0.$$
In particular, $W_{\emptyset}^{(p)}(\infty)$ converges to $W(\infty)$ in
$L^2(\P)$ as $p\rightarrow 1$.
\end{lemma}

In order to prove a joint limit law for
the processes $Z_i^{(p)}(\rho_n)$, $i\in\N$, we need information on their
birth times $b_i^{(n)}$ when $p\rightarrow 1$. This is achieved by the next
lemma, which corresponds to \cite[Lemma 4]{BeBr}.
\begin{lemma}
\label{L4BeBr}
As $p\rightarrow 1$, in the sense of finite-dimensional distributions,
$$
\left((1-p)W(\infty)\exp\left(pb_i^{(p)}\right) : i\in\N\right)\Longrightarrow \left(S_i :
i\in\N\right),$$
where $S_i = \ebf_1+\dots +\ebf_i$, and $\ebf_1, \ebf_2,\dots$ are
i.i.d. standard exponential variables.
\end{lemma}
As for Lemmas \ref{L2BeBr} and \ref{L3BeBr}, one can follow the proof
of \cite{BeBr} to get this last result for our model. As a consequence, we
see that for all $i\in\N$
\begin{equation}
\label{bi}
b_i^{(p)}= -(1/p)\ln(1-p) +O(1)\quad\hbox{as }p\rightarrow 1,
\end{equation}
i.e. the set of values $b_i^{(p)} + (1/p)\ln(1-p)$ is stochastically
bounded as $p\rightarrow 1$.

For the rest of this section, we let $p=p(n)$ depend on
$n$ such that $p\rightarrow 1$ as $n\rightarrow\infty$, but write mostly $p$
instead of $p(n)$. We first compute the asymptotic size
$Z_{\emptyset}^{(p)}(\rho_n)$ of the ancestral subpopulation, or, to put it
differently, the asymptotic size of the root cluster of a Bernoulli bond
percolation on $T_n$ with parameter $p$.
\begin{lemma}
\label{C2BeBr}
$$
\lim_{n\rightarrow\infty}n^{-p}Z_{\emptyset}^{(p)}(\rho_n)=1\quad\hbox{in probability.}
$$
\end{lemma}
\proof Since $\lim_{t\rightarrow\infty}\e^{-t}Z(t) = W(\infty)$ a.s., we
have $\lim_{n\rightarrow\infty}\e^{-\rho_n}n = W(\infty)$ a.s., implying
$$
\lim_{n\rightarrow\infty}\e^{-p\rho_n}
n^{p}W^{-1}(\infty)=1\quad\hbox{almost 
  surely.}
$$
On the other
hand, we know from Lemmas \ref{L2BeBr} and \ref{L3BeBr} that
$$
\lim_{n\rightarrow\infty}\e^{-\rho_n}Z(\rho_n) =
\lim_{n\rightarrow\infty}\e^{-p\rho_n}Z_\emptyset^{(p)}(\rho_n)=W(\infty)\quad\hbox{in probability.}
$$
This proves the statement.  \QED 

From the proof of the foregoing lemma, we deduce that {\it a fortiori}, $\rho_n = \ln
n + O(1)$ as $n\rightarrow\infty$. 

We recall that for Theorem \ref{T1} we
additionally require $(1-p(n))^{k}n^{p(n)}\rightarrow
\infty$ for all $k\in\mathbb{N}$. We will now implicitly assume
this at least for $k= 1$, so that in particular,
\begin{equation}
\label{tau-n}
\lim_{n\rightarrow\infty}(p\rho_n + \ln (1-p)) = \infty\quad\hbox{in probability}.
\end{equation}
We now consider two independent sequences $(W_i(\infty) : i\in\N)$ and $(\ebf_i
: i\in\N)$ of i.i.d. standard exponential random variables. We shall assume
that they are both defined on the same probability space. As before, let
$S_i = \ebf_1+\dots+\ebf_i$. Theorem $2$ of \cite{BeBr} tailored to
our needs yields finite-dimensional convergence of the $Z_i^{(p)}(\rho_n)$,
$i\in\N$.

\begin{proposition}
\label{T2BeBr}
As $n\rightarrow \infty$, in the sense of finite-dimensional 
distributions,
$$
\left(\frac{(1-p)^{-1}}{n^p}Z_i^{(p)}(\rho_n) :
  i\in\N\right)\Longrightarrow \left(\frac{W_i(\infty)}{S_i} : i\in\N\right).
$$
\end{proposition}
\proof
For $i\in\N$ and $t\geq 0$, put
$$ W_i^{(p)}(t) = \e^{-pt}Z_i^{(p)}(b_i^{(p)}+t).$$ Let 
$(t_i^{(p)})_{0<p<1}$ be a family of random times with $\lim_{p\rightarrow 1}t_i^{(p)}=\infty$ in
probability. From Lemma \ref{L3BeBr} (with $W_i^{(p)}$ and $W_i(\infty)$ in
place of $W_\emptyset^{(p)}$ and $W_\emptyset(\infty)$) we infer that there
is the convergence in the sense of finite-dimensional
distributions
$$
\left(Z_i^{(p)}(b_i^{(p)}+t_i^{(p)}) : i\in\N\right)\Longrightarrow
\left(W_i(\infty) : i\in\N\right).
$$
Concerning the birth times, we have by Lemma \ref{L4BeBr} the
finite-dimensional convergence
$$
\left(\frac{1}{(1-p)W(\infty)}\exp\left(-pb_i^{(p)}\right) : i\in\N\right)\Longrightarrow \left(\frac{1}{S_i} :
i\in\N\right).$$
By the remark below the definition of $b_u^{(p)}$, the sequence $(b_i^{(p)}:i\in\N)$ is independent from 
$(Z_i^{(p)}(b_i^{(p)}+\cdot): i\in\N)$, so that we have in fact joint weak 
convergence towards $(W_i(\infty), 1/S_i : i\in\N)$.
We finally let $t_i^{(p)}=\rho_n-b_i^{(p)}$. Then
$t_i^{(p)}\rightarrow\infty$ in probability for $p\rightarrow 1$ by
\eqref{bi} and \eqref{tau-n}. By the mapping theorem, the product of the
above left hand sides converges to $(W_i(\infty)/S_i : i\in\N)$, as claimed.
\QED

In order to obtain convergence of the tree of cluster sizes for the
first generation, we have to rank the sequence $(Z_i^{(p)}(\rho_n) :
i\in\N)$ in the decreasing order of their elements. Note that
finite-dimensional convergence for the reordered sequence cannot directly
be deduced from Proposition \ref{T2BeBr}. We first have to show that for
$\ell$ fixed, the $\ell$ largest subpopulations of generation $1$ at time
$\rho_n$ are with high probability to be found amongst the $k$ oldest when
$n\rightarrow\infty$ and $k\rightarrow\infty$. In view of the last
proposition and \eqref{bi}, we have to ensure that at time $\rho_n$, we see
only with small probability a subpopulation of size of order $(1-p)n^p$
which bears a single mutation and was born at a time much later than
$-(1/p)\ln(1-p)$.

For later use, namely for the proof of Corollary \ref{C1}, it will be
helpful to consider also subpopulations with more than one mutation. For
that purpose, let us list the full system of subpopulations
$(Z_u^{(p)}:u\in\pU)$ in the order of their birth times. We obtain a
sequence $(Y_i^{(p)}(t) : t\geq 0,\, i\in \N_0)$ which represents the same
process as $(Z_u^{(p)}(t): t\geq 0,\, u\in\pU)$, such that $Y_0^{(p)}=
Z_\emptyset^{(p)}$, $Y_1^{(p)}(b_1^{(p)}+\cdot)\equiv
Z_1^{(p)}(\cdot)$, and $Y_i^{(p)}(t)=0$ if less than $i$
mutants were born up to time $t$. Moreover, $(Z_i^{(p)}: i\in\N)$ is a
subsequence of $(Y_i^{(p)} : i\in\N_0)$ which corresponds to the
subpopulations with a single mutation. We denote by
$$
N^{(p)}(t) = \left|\left\{i\in\N: Y_i^{(p)}(t)>0\right\}\right| =
\left|\left\{u\in\pU\backslash\emptyset : b_u^{(p)}\leq t\right\}\right|
$$
the number of subpopulations born up to time $t$, discounting the ancestral
population of type $\emptyset$. Our next statement resembles \cite[Lemma
7]{BeBr}. However, in our setting we have to be more careful with the
estimates.
\begin{lemma}
\label{L7BeBr}
Let $\ve >0$. Then
$$
\lim_{r\rightarrow\infty}\lim_{n\rightarrow\infty}\P\left(\exists i\in\N :
  Y_i^{(p)}(-(1/p)\ln(1-p)+r) = 0\hbox{ and }Y_i^{(p)}(\rho_n)>\ve
  (1-p)n^p\right) = 0.
$$
\end{lemma}
\proof
Denote by $(\pF_t)_{t\geq 0}$ the natural filtration generated by the
system of processes $(Y_j^{(p)} : j\in\N_0)$. The counting process
$N^{(p)}$ is $(\pF_t)$-adapted, and its jump times $\gamma_k^{(p)} =
\inf\{t\geq 0: N^{(p)}(t) = k\}$ are $(\pF_t)$-stopping times.
By the strong Markov property, we see that each of the processes
$Y_k^{(p)}(\gamma_k^{(p)}+\cdot)$ for $k\in\N$ is a Yule process started
from a single particle of size $1$, with birth rate $p$ per unit population
size. Moreover, $Y_k^{(p)}(\gamma_k^{(p)}+\cdot)$ is independent of $\pF_{\gamma_k^{(p)}}$.
We let $r_n = -(1/p)\ln(1-p)+r$ and $s_n = \ln n +s$, where $r,s> 0$. The number of
processes $Y_k^{(p)}$ born after time $r_n$, which have at time $s_n$
a size greater than $\ve(1-p)n^p$ is given by
$$
X_n = \sum_{k=1}^\infty \1_{\left\{r_n< \gamma_k^{(p)}\leq
  s_n\right\}}\1_{\left\{Y_k^{(p)}(s_n)>\ve(1-p)n^p\right\}}=
  \int_{r_n}^{s_n}\1_{\left\{Y^{(p)}_{N^{(p)}(t)}(s_n)>\ve(1-p)n^p\right\}}\dt N^{(p)}(t).
$$
For $u\geq 0$ fixed, $Y_k^{(p)}(\gamma_k^{(p)}+u)$ is geometrically distributed
with parameter $\exp(-pu)$, see e.g. Yule \cite{Yu}. We obtain, using the bound
$(1-x)^a\leq \exp(-ax)$ in the second inequality,
\begin{align*}
  \E\left(X_n\right) &\leq \E\left(\int_{r_n}^{s_n}(1-\exp(-p(s_n-t)))^{\lfloor\ve(1-p)n^p\rfloor}\dt
    N^{(p)}(t)\right)\\
  &\leq
  \E\left(\int_{r_n}^{s_n}\exp\left[-(\ve/2)(1-p)n^p\exp(-p(s_n-t))\right]\dt
    N^{(p)}(t)\right).
\end{align*}
The dynamics of the family $(Y_i^{(p)}:i\in\N)$ and the strong Markov
property entail that $N^{(p)}$ grows at rate $(1-p)Z$, where
$Z=\sum_{u\in\pU}Z_u^{(p)}=\sum_{i=0}^\infty Y_{i}^{(p)}$ is a standard
Yule process. In particular, $N^{(p)}(t) - (1-p)\int_0^tZ(s)\dt s$ is a martingale.
Since $\E(Z(t)) = \e^t$, we get with the substitution
$x=\e^{pt}$ in the second line 
\begin{align*}
\E\left(X_n\right)&\leq(1-p)\int_{r_n}^{s_n}\e^t\exp\left[-(\ve/2)(1-p)n^p\exp(-p(s_n-t))\right]\dt
t\\
&=
\frac{1-p}{p}\int_{\e^{pr_n}}^{\e^{ps_n}}x^{(1-p)/p}\exp\left[-x(\ve/2)(1-p)n^p\exp(-ps_n)\right]\dt
x.
\end{align*}
We perform an integration by parts and substitute the values of $r_n$ and
$s_n$.  This gives
\begin{align*}
\E\left(X_n\right)&\leq
2\frac{\e^{ps}(1-p)^{(p-1)/p}}{\ve p}\e^{(1-p)r}\exp\left[-(\ve/2)
  \e^{p(r-s)}\right]\\
& +\;\;
2\frac{\e^{ps}(1-p)}{\ve p^2}\int_{\e^{pr_n}}^{\e^{ps_n}}x^{(1-2p)/p}\exp\left[-x(\ve/2)(1-p)\exp(-ps)\right]\dt x.
\end{align*}
For large $n$, we have $p>1/2$ and $x^{(1-2p)/p} \leq 1$ on the domain of
integration. For such $n$
$$
\E\left(X_n\right)\leq
(4/\ve)(1-p)^{(p-1)/p}\e^{(1-p)r+ps}\exp\left[-(\ve/2)
  \e^{p(r-s)}\right]+(16/\ve^2)\e^{2ps}\exp\left[-(\ve/2)
  \e^{p(r-s)}\right].
$$
In particular, for fixed $r,s,\ve>0$,
$$
\limsup_{n\rightarrow\infty}\E\left(X_n\right)\leq 16(\ve^{-1} +\ve^{-2})\e^{2s}\exp\left[-(\ve/2)
  \e^{(r-s)}\right],
$$
and the right side converges to zero when $r\rightarrow\infty$ and $\ve$,
$s$ are fixed. We have shown that
$$
\lim_{r\rightarrow\infty}\lim_{n\rightarrow\infty}\P\left(\exists i\in\N :
  Y_i^{(p)}(-(1/p)\ln(1-p)+r) = 0\hbox{ and }Y_i^{(p)}(s_n)>\ve
  (1-p)n^p\right) = 0.
$$
Since $\lim_{s\rightarrow\infty}\P\left(\rho_n>\ln n +s\right) = 0$ by
\eqref{tau-n}, the lemma is proved.  
\QED

We are now in position to prove Theorem \ref{T1} restricted to
generations $0$ and $1$. We recall that $0<p=p(n)<1$ with $p\rightarrow 1$ and
$(1-p)n^p\rightarrow\infty$ as $n\rightarrow\infty$.

\begin{proposition}
\label{Pgen01conv}
As $n\rightarrow\infty$, $\pC_{\emptyset}^{(n)}\sim n^p$ in
probability, and for every fixed $\ell\in\N$,
$$
\left(\frac{(1-p)^{-1}}{n^{p}}\pC_1^{(n)},\dots,\frac{(1-p)^{-1}}{n^{p}}\pC_{\ell}^{(n)}\right)\Longrightarrow
(\pZ_1,\dots,\pZ_{\ell}).
$$
\end{proposition}

\proof The convergence of the root cluster was already shown in Lemma
\ref{C2BeBr}. Indeed, it follows from our construction that
$Z_{\emptyset}^{(p)}(\rho_n)$ is distributed as the size
$\pC_\emptyset^{(n)}$ of the root cluster of a Bernoulli bond percolation
on $T_n$ with parameter $p(n)$. Next, we deduce from Proposition
\ref{T2BeBr} and \eqref{bi} together with Lemma \ref{L7BeBr} that if we
write $(x_i)^\downarrow$ for the decreasing rearrangement of a sequence of
positive real numbers $(x_i)$ with either pairwise distinct elements or
finitely many non-zero terms, we have
\begin{equation}
\label{gen1limit}
\left(\frac{(1-p)^{-1}}{n^p}Z_i^{(p)}(\rho_n) :
  i\in\N\right)^\downarrow\Longrightarrow \left(\frac{W_i(\infty)}{S_i} : i\in\N\right)^\downarrow
\end{equation}
in the sense of finite-dimensional distributions as $n$ tends to
infinity. Since
$$
\left(Z_i^{(p)}(\rho_n) :
  i\in\N\right)^\downarrow \eql \left(\pC_i^{(n)} : i\in\N\right),
$$
it only remains to identify the limit on the right hand side of
\eqref{gen1limit}. Recalling that $(S_i:i\in \N)$ is independent of
$(W_i(\infty) : i\in \N)$, $((S_i,W_i(\infty)): i\in \N)$ can be viewed as
the sequence of atoms of a Poisson point process on
$(0,\infty)\times(0,\infty)$ with intensity $\dt s\otimes \e^{-r}\dt r$, and
the claim follows.\QED

We now turn to the proof of Corollary \ref{C1} stated in Section
\ref{Smainresults}. In view of what we have already proved, it will be
sufficient to check that for each fixed $r>0$, the subpopulations which
were born up to time $-(1/p)\ln(1-p)+r$ carry all a single mutation with
high probability when $n\rightarrow\infty.$ Similarly to the definition of $N^{(p)}(t)$, we let 
$$M^{(p)}(t) = \left|\left\{i\in\N : b_i^{(p)}\leq t\right\}\right|$$
denote the number of subpopulations with a single mutation at time $t$. The
following statement is similar to Lemma $6$ in \cite{BeBr}.
\begin{lemma}
\label{L6BeBr}
Let $\Delta^{(p)}(t) = N^{(p)}(t)-M^{(p)}(t)\geq
0$ denote the number of subpopulations born up to time $t$, which bear more than one single
mutation. Then for each $r>0$,
$$
\lim_{n\rightarrow\infty}\E\left(\Delta^{(p)}(-(1/p)\ln(1-p)+r)\right)=0.
$$
\end{lemma}

\proof
Let $r_n = -(1/p)\ln(1-p)+r$. Since $N^{(p)}(t) -
(1-p)\int_0^tZ(s)\dt s$ is a martingale, 
$$
\E\left(N^{(p)}(r_n)\right) = (1-p)\int_0^{r_n}\E(Z(s))\dt s.
$$
Similarly, we obtain
$$
\E\left(M^{(p)}(r_n)\right) = (1-p)\int_0^{r_n}\E\left(Y_0^{(p)}(s)\right)\dt s.
$$
Using $\E(Z(s))= \e^s$, $\E(Y_0^{(p)}(s))=\e^{ps}$ and $p(n)\rightarrow 1$,
a small computation shows $\E\left(\Delta^{(p)}(r_n)\right)=o(1)$ for
$n\rightarrow\infty$.  \QED

\proofof{\bf Corollary \ref{C1}:}\hskip10pt From \eqref{bi}, Lemmas
\ref{L7BeBr}, \ref{L6BeBr} and Proposition \ref{T2BeBr}, we see that the
sizes of the largest non-ancestral subpopulations at time $\rho_n$ are
taken by the subpopulations with a single mutation only. Recalling the
connection between subpopulations at time $\rho_n$ and percolation
clusters, the proof of the corollary is then a consequence of Proposition
\ref{Pgen01conv}. \QED

\subsection{Higher generation convergence}
\label{SproofT1}
The recursive structure of $T_n$ allows us to transfer the arguments
of the foregoing section to higher generation clusters. We however need some
preparation.

Let $T_n$ be a RRT on $\{0,1,\dots,n\}$, as usual. Here it will be convenient 
to label the edges of $T_n$ by their outer endpoints, i.e. the edge $e$
joining vertex $i$ to vertex $j$, where $i<j$, is labeled $j$. We then say
that $e$ is the $j$th edge of $T_n$.  We incorporate Bernoulli bond
percolation on $T_n$, but instead of deleting edges, we simply mark them
with probability $1-p$ each, independently of each other. After such a
marking of edges, we call a subtree of $T_n$ {\it intact}, if it contains
only unmarked edges and is maximal in the sense that no further
edges without marks can be attached to it. In other words, the intact
subtrees of $T_n$ are precisely the percolation clusters of $T_n$.

We again view $T_n$ as the genealogical tree of a standard Yule process
stopped at the instant when the $(n+1)$th individual is born. Henceforth we
will identify vertices with individuals, i.e. we will make no difference
between the vertex labeled $j$ and the $j$th individual of the population
system. The marked edges indicate a birth event of a mutant.  This means
that if the $j$th edge is a marked edge, then the $j$th individual is a
mutant, and the vertices of the intact subtree rooted at $j$ correspond to
the individuals bearing the same genetic type as the $j$th
individual. Moreover, the genetic type $u\in\pU$ of the $j$th individual can be derived 
from the subtree of $T_n$ spanned by the vertices $0,1,\dots,j$ and from the marks on
its edges.

Our description shows that we may generate the subpopulation
sizes $Z_u^{(p)}(\rho_n)$, $u\in \pU$, by first 
picking a RRT $T_n$, then marking each edge with probability $1-p$,
independently of each other, and then defining $Z_u^{(p)}(\rho_n)$ to be the size
of the intact subtree of $T_n$ rooted at the mutant of type $u$.

Let us write $\tau_u^{(n)}$ for the full genealogical (sub)tree which stems
from the mutant of type $u$. This means that $\tau_u^{(n)}$ is the maximal
subtree of $T_n$ rooted at the mutant of type $u$, including all marked and
unmarked edges above its root. Clearly, $\tau_u^{(n)}$ might contain
several intact subtrees of $T_n$, and we agree that
$\tau_{\emptyset}^{(n)}$ is given by $T_n$ itself. Moreover, we let
$\tau_u^{(n)}=\emptyset$ if there is no mutant of type $u$. For example,
the non-empty vertex sets of the genealogical subtrees of the recursive
tree on the left side of Figure $1$ (dashed lines represent the marked
edges) are given by $\tau_\emptyset = \{0,1,\dots,10\}$, $\tau_1=\{2\}$,
$\tau_2=\{3,5,7,8,9,10\}$, $\tau_{21}=\{5,9\}$, $\tau_{22}=\{10\}$.

Let us introduce the following terminology. For an arbitrary subset
$A\subset\{0,1,\dots,n\}$ of size $k$, we call the bijective map from $A$
to $\{0,1, \dots,k-1\}$, which preserves the order, the {\it canonical
relabeling} of vertices. Clearly, the canonical relabeling transforms a
recursive tree on $A$ into a recursive tree on $\{0,1,\dots, k-1\}$.

We next observe that conditionally on its size $|\tau^{(n)}_u|=k$ and upon
the canonical relabeling of its vertices, $\tau^{(n)}_u$ is itself
distributed as a RRT on $\{0,1,\dots, k-1\}$. Indeed, as we pointed out
above, in order to decide whether a given vertex $j$ of $T_n$ is the root
of the subtree encoded by $\tau^{(n)}_u$, we have to look only at the
subtree (with its marks) spanned by the vertices $\{0,1,\dots,j\}$. In
particular, the structure of the subtree stemming from $j$ is
irrelevant. If we condition on $|\tau^{(n)}_u|=k$ and perform the canonical
relabeling of vertices, the recursive construction then implies that each
increasing arrangement of the vertices $\{0,1,\dots,k-1\}$ is equally
likely, that is to say $\tau^{(n)}_u$ is a random recursive tree. Moreover,
if $u,v\in\pU$ do not lie on the same infinite branch of $\pU$ emerging
from the root $\emptyset$, then $\tau^{(n)}_u$ and $\tau^{(n)}_v$ are
conditionally on their sizes independent RRT's, since their vertex sets are
disjoint.

We remark that these properties of random recursive trees are closely
related to the so-called {\it splitting property}, which plays a major role
in our analysis of the destruction process of a RRT in Section
\ref{Streeofcomponents}.

We henceforth call a subtree $\tau_u^{(n)}$ with $|u|=k$ a subtree of
generation $k$. Our final main step for proving Theorem \ref{T1} is a
convergence result for the ``tree of subtree sizes'' $(|\tau_u^{(n)}| :u\in
\pU)$. For that purpose, we decreasingly order the children of each element
$|\tau_u^{(n)}|$, but keep the parent-child relation. More precisely, each
element $|\tau_u^{(n)}|$ has finitely many non-zero children, say
$|\tau_{u1}^{(n)}|,\dots,|\tau_{u\ell}^{(n)}|$, and we let
$\sigma_u:\N\rightarrow\N$ be the random bijection which sorts this
sequence in the decreasing order, i.e.
\begin{equation*}
  |\tau_{u\sigma_u(1)}^{(n)}|\geq |\tau_{u\sigma_u(2)}^{(n)}|\geq\dots\geq|\tau_{u\sigma_u(\ell)}^{(n)}|,
\end{equation*}
with $\sigma_u(i) = i$ for $i>\ell$.

Out of these maps we define the global random bijection $\sigma =
\sigma^{(n)} : \pU\rightarrow\pU$ recursively by setting $\sigma(\emptyset)
= \emptyset$, $\sigma(j)=\sigma_\emptyset(j)$, and then, given $\sigma(u)$,
$\sigma(uj) = \sigma(u)\sigma_{\sigma(u)}(j)$, $u\in\pU$,
$j\in\mathbb{N}$. Note that $\sigma$ indeed preserves the parent-child relation,
i.e. children of $u$ are mapped into children of $\sigma(u)$. 

We recall that $p=p(n)\rightarrow 1$, and
$(1-p(n))^{k}n^{p(n)}\rightarrow\infty$ for each $k\in\N$.

\begin{proposition}
\label{Ptreeofsubtreesizes}
As $n\rightarrow\infty$, in the sense of
finite-dimensional distributions,
$$
\left(\frac{(1-p(n))^{-|u|}}{n}|\tau_{\sigma(u)}^{(n)}| :u\in\pU\right)\Longrightarrow
(\pZ_u : u\in\pU).$$
\end{proposition}
\proof The convergence of $|\tau_{\emptyset}^{(n)}|/n$ is trivial. Let us
first show that in the sense of finite-dimensional distributions, as
$n\rightarrow\infty$,
\begin{equation}
\label{conv1}
\left(\frac{(1-p)^{-1}}{n}|\tau_{\sigma(i)}^{(n)}| : i\in\N\right)\Longrightarrow
\left(\pZ_i : i\in\N\right).
\end{equation}
Fix $\ell\in\N$ and denote by $\tilde{\pC}_i^{(n)}$ the size of the root
percolation cluster inside $\tau_i^{(n)}$, i.e. the size of the intact
subtree with the same root node as $\tau_i^{(n)}$. Since conditionally on
its size (and upon the canonical relabeling), $\tau_i^{(n)}$ is a random recursive tree, Lemma \ref{C2BeBr}
shows that for each $i=1,\dots,\ell$, $$|\tau_{i}^{(n)}|\sim
\left(\tilde{\pC_{i}}^{(n)}\right)^{1/p}\quad\hbox{in probability.}$$
Furthermore, in the notation from Section \ref{SYule} we have the equality
in distribution
$$
\left(\tilde{\pC}_1^{(n)},\dots,\tilde{\pC}_{\ell}^{(n)}\right)\eql
\left(Z_1^{(p)}(\rho_n),\dots,Z_\ell^{(p)}(\rho_n)\right).
$$
Proposition \ref{T2BeBr} (with $W_i(\infty)$,
$S_i$ defined there), the last two displays and the fact that $p\rightarrow
1$ imply the convergence in distribution 
$$
\left(\frac{(1-p)^{-1}}{n}|\tau_{1}^{(n)}|,\dots,\frac{(1-p)^{-1}}{n}|\tau_{\ell}^{(n)}|\right)\Longrightarrow
\left(\frac{W_1(\infty)}{S_1},\dots,\frac{W_\ell(\infty)}{S_\ell}\right).
$$
From Lemma \ref{C2BeBr} and \eqref{bi} together with Lemma
\ref{L7BeBr}, we know that the $\ell$ largest subtrees amongst
$\tau_i^{(n)}$, $i\in\N$, are with high probability to be found under the
first $k$, provided $k$ and $n$ are large. With our identification of the ranked
sequence $(W_i(\infty)/S_i : i\in\N)$ from the proof of Proposition
\ref{T2BeBr}, the last display therefore implies \eqref{conv1}.

We next show that in the sense of finite-dimensional distributions,
\begin{equation}
\label{conv2}
\left(\frac{(1-p)^{-|u|}}{n}|\tau_{\sigma(u)}^{(n)}| :u\in\pU,|u|\leq
  2\right)\Longrightarrow \left(\pZ_u :u\in\pU,|u|\leq 2\right).
\end{equation} 
As we already remarked, for disjoint integers $j,k$, the RRT's
$\tau_j^{(n)}$ and $\tau_k^{(n)}$ are conditionally on their sizes {\it
  independent} RRT's. Since we have just proved finite-dimensional
convergence of $(|\tau_{\sigma(j)}^{(n)}|:j\in\N)$, it will therefore be
enough to show that for $g,f^{(i)}:[0,\infty)\rightarrow [0,1]$ bounded and
uniformly continuous, and $j,\ell \in\N$,
\begin{align*}
  \lefteqn{\E\left[g\left(\frac{(1-p)^{-1}}{n}|\tau_{\sigma(j)}^{(n)}|\right)f^{(1)}\left(\frac{(1-p)^{-2}}{n}|\tau_{\sigma(j1)}^{(n)}|\right)\dots
      f^{(\ell)}\left(\frac{(1-p)^{-2}}{n}|\tau_{\sigma(j\ell)}^{(n)}|\right)\right]}\\
  &\rightarrow \E\left[g(\pZ_j)f^{(1)}(\pZ_ja_{1})\dots
    f^{(\ell)}(\pZ_ja_{\ell})\right],\hspace{3cm}
\end{align*} 
where for $i=1,\dots,\ell$, $a_{i}$ is the $i$th
largest atom of a Poisson random measure on $(0,\infty)$ with intensity
$\nu(da)=a^{-2}\dt a$. For ease of readability, we restrict ourselves to
the case $\ell=1$, the case $\ell\geq 2$ being similar. By the properties
of $\tau_u^{(n)}$ discussed above, we have for each $m=0,\dots,n-j$,
\begin{equation*}
 \E\left[f^{(1)}\left(\frac{(1-p)^{-2}}{n}|\tau_{\sigma(j1)}^{(n)}|\right)\Big|\, |\tau_{\sigma(j)}^{(n)}|=m+1\right] =
  \E_m\left[f^{(1)}\left(\frac{(1-p)^{-2}}{n}|\tau_{\ast}^{(m)}|\right)\right], 
\end{equation*}
where $\E_m$ is the mathematical expectation starting from a random
recursive tree $T_m$ with $m+1$ vertices, and under $\E_m$,
$\tau_{\ast}^{(m)}$ is the largest amongst all (full) genealogical
subtrees of $T_m$ stemming from a mutant with one single mutation.

Now if $m\sim (1-p)na$ for some fixed $a>0$, $m=m(n)$ integer-valued, we
obtain from the convergence of generation $1$,
i.e. \eqref{conv1}, that
\begin{equation*}
 \E_m\left[f^{(1)}\left(\frac{(1-p)^{-2}}{n}|\tau_{\ast}^{(m)}|\right)\right] \sim \E\left[f^{(1)}(aa_1)\right],
\end{equation*} 
with $a_1$ the largest atom of a Poisson random measure with intensity
$\nu$.  Since we have already proved that
$(1-p)^{-1}n^{-1}|\tau_{\sigma(j)}^{(n)}|$ converges in distribution to $\pZ_j$,
and since the map $ (a,a_1)\mapsto g(a)f^{(1)}(aa_1)$ is uniformly
continuous on bounded sets, this establishes~\eqref{conv2}. With the same
arguments, we obtain finite-dimensional convergence up to generation $3$,
then up to $4$, and so on, so that the proposition is proved.\QED

We finish now the proof of Theorem \ref{T1}.

\proofof{\bf Theorem \ref{T1}:}\hskip10pt 
Write $\tilde{\pC}_u^{(n)}$ for the size of the root percolation cluster
inside $\tau_u^{(n)}$. From Lemma \ref{C2BeBr} we know that as $n\rightarrow\infty$,
$\tilde{\pC}_u^{(n)}\sim |\tau_u^{(n)}|^p$ in probability. Since 
$(1-p)^{-kp}\sim(1-p)^{-k}$ as $n\rightarrow\infty$ for each $k\in\N$, we obtain from Proposition
\ref{Ptreeofsubtreesizes}
$$
\left(\frac{(1-p)^{-|u|}}{n^p}\tilde{\pC}_{\sigma(u)}^{(n)} :
    u\in\pU\right)\Longrightarrow \left(\pZ_u:u\in\pU\right)\quad\hbox{as }n\rightarrow\infty
$$
in the sense of finite-dimensional distributions. Now observe that the
process $(\tilde{\pC}_{\sigma(u)}^{(n)} : u\in\pU)$ does already encode all
cluster sizes and the genealogical structure of the corresponding clusters.
It is however possible that $\tilde{\pC}_{\sigma(ui)}^{(n)}<
\tilde{\pC}_{\sigma(uj)}^{(n)}$ for $i<j$, since the map $\sigma$ provides
an ordering according to the sizes of the surrounding subtrees
$\tau_u^{(n)}$, not according to the cluster sizes. 

Therefore, in order to finish the proof, we need to argue that the
convergence remains true if the family $(\tilde{\pC}_{\sigma(u)}^{(n)} :
u\in\pU)$ is ranked in the decreasing order, that is if the children
$(\tilde{\pC}_{\sigma(ui)}^{(n)}:i\in\N)$ of each element
$\tilde{\pC}_{\sigma(u)}^{(n)}$ are decreasingly sorted according to their
sizes, under preservation of the parent-child relation. However, Proposition
\ref{Ptreeofsubtreesizes} and the fact that
$\tilde{\pC}_{u}^{(n)}\sim |\tau_{u}^{(n)}|^{p}$ in probability 
entail that we find the $\ell$ largest elements of
$(\tilde{\pC}_{\sigma(ui)}^{(n)}:i\in\N)$ amongst the first $k$ with 
probability as close to $1$ as we wish, provided we choose $k$ and $n$ large enough.
This completes the proof of the theorem.
\QED

The tree $\pC^{(n)}$ provides arguably the most natural encoding of the cluster
sizes. From the point of view of a dynamical version of a Bernoulli bond
percolation, where the edges are removed one after the other in a random
uniform order, another possibility might however come to ones mind, which
takes into account the order in which the edges are removed. We will 
discuss this in Section \ref{Sapps} for the supercritical regime. Our
discussion is based on a result for the so-called {\it tree of
  components}, which we present next.

\section{The tree of components}
\label{Streeofcomponents}
\subsection{Our main result on the destruction of a RRT}
\label{Sdestruction}
We consider on $T_n$ a continuous-time destruction process with parameter
$1/\ln n$. This means that we attach to each edge
$e$ of $T_n$ an independent exponential clock $\ebf(e)$ of parameter
$1/\ln n$, and we delete edge $e$ at time $\ebf(e)$. After the $n$th edge
has been deleted, the tree has been completely
destructed, and the process terminates.

We encode the sizes and birth times of the tree components stemming from
the destruction of $T_n$ by a tree-indexed process
$(\pB^{(n)},b^{(n)})=((\pB^{(n)}_u,b^{(n)}_u) : u\in\pU)$, the {\it tree of
  components}. Here, $\pU$ denotes again the universal tree.

As it is already apparent from \cite{Be1}, choosing the parameter in the
destruction process equals $1/\ln n$ is the most natural choice, since it
leads without normalization to a simple description of the birth times in the limit
$n\rightarrow\infty$.

Following the steps of the destruction process, we build the tree of
components dynamically starting from the singleton $(\pB^{(n)}_\emptyset,
b^{(n)}_\emptyset)=(n+1, 0)$ and ending after the $n$th edge removal with
the full process $(\pB^{(n)},b^{(n)})$. For the ease of readability, we
omit the superscript $(n)$ in the following construction.

Let $e_1,\dots,e_n$ denote the edges of $T_n$ listed in
the increasing order of their exponential clocks, i.e. such that
$\ebf(e_1)<\ebf(e_2)<\dots<\ebf(e_n)$. Then at time $\ebf(e_1)$,
$e_1$ is the first edge to be removed from $T_n$, and $T_n$ splits into two
subtrees, say $\tau_n^0$ and $\tau_n^\ast$, where $\tau_n^0$ contains the
root $0$. The size $|\tau_n^\ast|$ is viewed
as the first child of $\pB_\emptyset$ and denoted by $\pB_1$,
and $b_1$ is set to $\ebf(e_1)$. Now first suppose that $e_2$ connects
two vertices in $\tau_n^\ast$. Then, at time $\ebf(e_2)$, $\tau_n^\ast$
splits into two tree components. The size of the component not containing
the root of $\tau_n^\ast$ is viewed as the first child of $\pB_1$ and
denoted by $\pB_{11}$, and $b_{11}$ is set to $\ebf(e_2)$. On
the contrary, if $e_2$ connects two vertices in $\tau_n^0$, then the size
of the component not containing $0$ is viewed as the second child of
$\pB_\emptyset$ and denoted by $\pB_2$, while $b_2$ is
set to $\ebf(e_2)$. 

It should be plain how to iterate this construction. After the $n$th edge
removal, we have in this way defined $n+1$ pairs of variables $(\pB_u,
b_u)$ with $|u|\leq n$. We extend the definition to the full universal tree
by letting $(\pB_u, b_u)=(0,\infty)$ for all the remaining $u\in\pU$.  Note
that for a non-zero element $\pB_u$, we have
$\pB_{ui}<\pB_{u}$ and $b_u < b_{ui}$ for $i\in\N$ (children are strictly smaller than
their parent, and they cannot be born before their parent), and $b_{ui}<
b_{uj}$ if $i<j$ and $b_{ui}\neq \infty$.

We point out that $(\pB,b)$ represents a final state of the destruction
when all edges have been deleted. In particular, the tree of components is
a process indexed by $\pU$, not by the time.  We call the first-coordinate
process $\pB=(\pB_u:u\in\pU)$ simply the {\it tree of component sizes}. An
example is shown in Figure $2$. Note that unlike our convention from
Section \ref{SproofT1} of labeling the edges by their outer endpoints, the
edge labels are here given by the order in which the edges were removed.
Similar to above, we say that a tree component is a component of generation
$k$ if its size is encoded by an element $\pB_u$ with $|u|=k$.
\begin{figure}[ht]
\begin{center}\parbox{7cm}{\includegraphics[width=6cm]{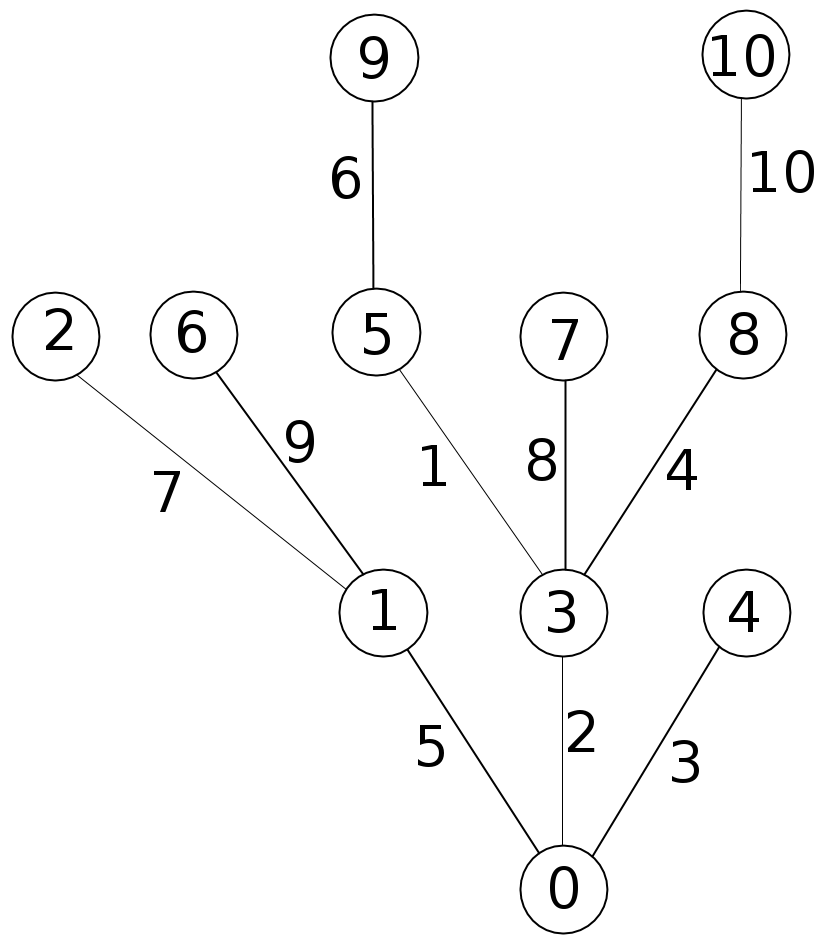}}
\parbox{9cm}{\includegraphics[width=9cm]{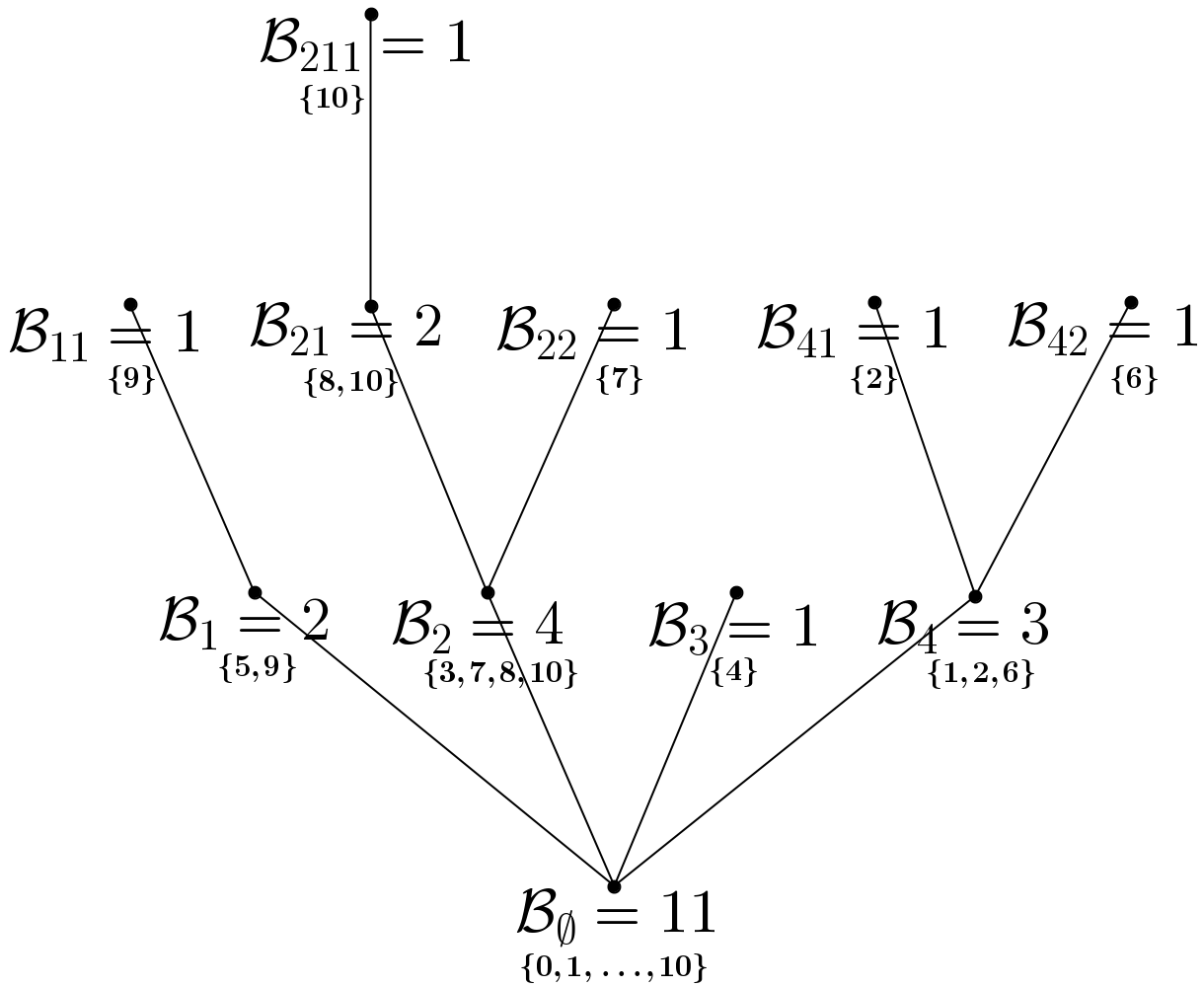}}
\end{center}
\centerline{\bf Figure 2} \centering{\sl Left: A recursive tree with
  vertices labeled $0,1,\dots,10$. The labels on the edges indicate the
  order in which they
  were removed by the destruction process.}\\
{\sl Right: The corresponding tree of component sizes, with the vertex sets
  of the tree components. The elements $\pB_u$ of size $0$ are omitted.}
\end{figure}

For our limit result, we will rank the children of each element
$(\pB_u^{(n)},b_u^{(n)})$ in the decreasing order of their {\it first}
coordinate, in the same way as we ranked the subtree sizes in the foregoing section. Say
$\pB_{u1}^{(n)},\dots,\pB_{u\ell}^{(n)}$ are the non-zero children of
$\pB_u^{(n)}$. We choose a random bijection $\sigma_u:\N\rightarrow\N$ with
$\sigma_u(i) = i$ for $i>\ell$ such that
\begin{equation*}
  \pB^{(n)}_{u\sigma_u(1)}\geq \pB^{(n)}_{u\sigma_u(2)}\geq \dots\geq\pB^{(n)}_{u\sigma_u(\ell)}.
\end{equation*}
As before, $\sigma = \sigma^{(n)} : \pU\rightarrow\pU$ is defined recursively by
setting $\sigma(\emptyset) = \emptyset$, $\sigma(j)=\sigma_\emptyset(j)$,
and given $\sigma(u)$, $\sigma(uj) = \sigma(u)\sigma_{\sigma(u)}(j)$,
$u\in\pU$, $j\in\mathbb{N}$. 

The limit object of the tree of component sizes is again given by the 
tree-indexed process $(\pZ_u:u\in\pU)$ from Section
\ref{Streeofclustersizes}. The birth times of the corresponding tree
components are in the limit described by a sequence $(\ebf_u :
u\in\pU\backslash\{\emptyset\})$ of i.i.d. standard exponential random
variables, which are independent of $\pZ$. We let $z_\emptyset =
\ebf_\emptyset = 0$, and, writing $u=(u_1\dots u_k)$ for $u$ of length
$k\geq 1$,
\begin{equation*}
z_u = \ebf_{u_1} + \ebf_{u_1u_2} +\dots + \ebf_{u_1\dots u_k}.
\end{equation*}
In words, $z_u$ is the sum of all the $\ebf$-values along the path in $\pU$ from
the root $\emptyset$ to $u$. 

We shall prove the following finite-dimensional convergence.
\begin{theorem}\label{T2} 
As $n\rightarrow\infty$, there is the
convergence in the sense of finite-dimensional distributions,
$$
\left(\left(\frac{(\ln n)^{|u|}}{n}\pB_{\sigma(u)}^{(n)},
    b_{\sigma(u)}^{(n)}\right) : u\in\pU\right)\Longrightarrow
  \left(\left(\pZ_u, z_u\right) : u\in\pU\right).
$$
\end{theorem}
For our application to (supercritical) Bernoulli bond percolation on $T_n$, it will be useful to
consider only those tree components which are born before a certain finite
time $t(n)\rightarrow t>0$. In this case, the limit object is 
obtained from ``squeezing-out'' the elements $(\pZ_u, z_u)$ with $z_u\geq t$. The 
corresponding limit statement is Proposition \ref{P2}, which is deferred to Section \ref{Sapps}.

\subsection{Some features of RRT's and of the tree of components}
\label{Stools}
\subsubsection{The splitting property}
As for the higher generation convergence of cluster sizes in Section
\ref{SproofT1}, we use the fractal structure of random recursive trees to
study the destruction process. Here we shall employ the splitting property,
which informally states that if an edge is removed uniformly at random from
a RRT, then the two subtrees which were connected by this edge are,
conditionally on their sizes, independent RRT's.

More precisely, select an edge of $T_n$ uniformly at random and remove
it. Then $T_n$ splits into two subtrees, say $\tau_n^0$ and $\tau_n^\ast$,
where $\tau_n^0$ contains the root.  Let $\xi$ be the integer-valued
variable with distribution
\begin{equation}\label{Estepdistr}
  \P(\xi=j)=\frac{1}{j(j+1)},\quad j=1,2, \dots.
\end{equation}
Remember that we call the canonical relabeling of a vertex set
$A\subset\{0,1,\dots,n\}$ the relabeling of its vertices by the labels
$0,1,\dots,|A|-1$, according to the increasing order of the original labels.
\begin{lemma}{\rm (Meir and Moon \cite{MM2})}
  \label{Lsplittingprop}
  Conditionally on the size $|\tau^0_n|=k$, the subtrees $\tau^0_n$ and
  $\tau^{\ast}_n$ are upon the canonical relabeling of their vertices
  independent random recursive trees on $\{0,1,\dots,k-1\}$ and
  $\{0,1,\dots,n-k\}$, respectively. Moreover, in the notation from above,
  $|\tau_n^\ast|$ has the same law as $\xi$ conditioned on $\xi\leq n$,
  that is
  \begin{equation}\label{Esplittingproba}\P(|\tau^{\ast}_n|=j)=\frac{n+1}{nj(j+1)},\quad
    j=1, \dots, n.\end{equation}

\end{lemma}
The first property in the statement is also referred to as the {\it
  splitting property} of random recursive trees.

\subsubsection{A coupling of Iksanov and M\"ohle}
Meir and Moon \cite{MM2} proved Equation \eqref{Esplittingproba} in their
study of the number $X_n$ of random cuts that are needed to isolate the
root $0$ in $T_n$. They considered the following algorithm for isolating
vertex $0$. Start from $T_n$ and remove an edge chosen uniformly at
random. Then iterate the procedure with the subtree which contains $0$, and
so on, until after $X_n\leq n$ steps, the root is finally isolated.

In \cite{MM2} Meir and Moon obtained first and second moment estimates for
$X_n$ and showed that $\lim_{n\rightarrow\infty}(\ln n/n)X_n=1$ in
probability. Later, Drmota {\it et al.} \cite{DIMR} proved a weak
limit law, showing that $n^{-1}(\ln n)^2X_n-\ln n-\ln\ln n$ converges in
distribution to a completely asymmetric Cauchy variable. A short
probabilistic proof of this result was found by Iksanov and M\"ohle \cite{IM}, which
turns out to be particular useful for our purpose. It is based on a
coupling of component sizes with the steps of an increasing random
walk. More precisely, let $\xi_1, \xi_2, \dots$ be a sequence of
i.i.d. copies of $\xi$, see \eqref{Estepdistr}, and set
$S_0=0$,
\begin{equation}\label{ESn}
S_n=\xi_1+\dots + \xi_n.
\end{equation} 
Denote the last time the random walk $S$ remains below level $n$ by
$$
L(n)=\max\{k\geq 0: S_k\leq n\}.
$$

\begin{lemma} \label{L2}{\rm (Iksanov and M\"ohle \cite{IM})} One can
  construct on the same probability space a random recursive tree $T_n$
  together with the random algorithm for isolating the root, and a version
  of the random walk $S$, such that if
$$
\pB_1^{(n)},\pB_2^{(n)},\dots,\pB_{X_n}^{(n)}
$$
denotes the sequence of the sizes of the subtrees which are cut off from
 the root component one after the other by the algorithm, then
 $X_n\geq L(n)$ and
$$
\left(\pB_1^{(n)},\dots, \pB_{L(n)}^{(n)}\right) = \left(\xi_1,\dots, \xi_{L(n)}\right).
$$
\end{lemma}
Besides the coupling, we use the following two facts about the random walk $S$ and its
last passage time $L(n)$, which can be found in a stronger form in
\cite{IM}.
\begin{equation}
\label{EQLn}
\lim_{n\rightarrow\infty}\frac{\ln n}{n}L(n)=1\quad\hbox{and}\quad
\lim_{n\rightarrow\infty}\frac{\ln n}{n}(n-S_{L(n)})=0\quad\hbox{in
  probability}.
\end{equation}
Combined with extreme value theory, we will use the coupling to determine the
asymptotic sizes and birth times of the tree components.

\subsubsection{Branching property of the tree of components}
The tree-indexed process $\pB^{(n)}$ can be interpreted as the genealogical
tree of a multi-type population model, where the type reflects the size of
the tree component. In particular, the ancestor $\emptyset$ has type
$n+1$. We stress that the characteristic ``type'' is used here in a
different way than in Section \ref{SYule}. A node $u$ with $\pB^{(n)}_u=0$
corresponds to an empty component and is therefore absent in the population
model. The splitting property leads to the following description.
\begin{lemma}
\label{L3}  
The population model induced by the tree of component sizes $\pB^{(n)}$ is
a multi-type Galton-Watson process starting from one particle of type
$n+1$. The reproduction distribution $\lambda_i$ of an individual of type
$i\geq 1$ is given by the law of the sequence of the sizes of the non-root
subtrees which are produced in the algorithm for isolating the root of a
RRT of size $i$.
\end{lemma}
We remark that the type of an individual is simply given by the total
size of the subtree of the genealogical tree stemming from that
individual. Therefore, types can be recovered from the sole structure of
the genealogical tree. 

When we incorporate the birth times of the tree components, Lemma
\ref{L3} and basic properties of exponential variables immediately yield
the following branching property.

\begin{lemma}\label{L4}
For every integer $k\geq 0$, conditionally on $((\pB_u^{(n)}, b_u^{(n)}) : |u|\leq
k)$, the families of variables
\begin{equation*}
\left(\left(\pB^{(n)}_{uj}, b_{uj}^{(n)}\right): j\in\mathbb{N}\right),\quad u\in\pU \mbox{ of length }|u|=k,
\end{equation*}
are independent, and the conditional law of each family
$((\pB^{(n)}_{uj}, b_{uj}^{(n)}): j\in\mathbb{N})$ only
depends on $(\pB^{(n)}_{u}, b_{u}^{(n)})$. 
More precisely, given $(\pB^{(n)}_u, b^{(n)}_{u})$ with
$\pB^{(n)}_{u}-1=m\geq 2$, there is the equality in
distribution
$$\left(\left(\pB^{(n)}_{uj},
    b^{(n)}_{uj}\right):j\in\mathbb{N}\right)\eql
\left(\left(\pB^{(m)}_j, b_u^{(n)}+\frac{\ln n}{\ln
      m}b_j^{(m)}\right):j\in\mathbb{N}\right),$$ where $\pB^{(m)}_j$ and
$b_j^{(m)}$ stem from a destruction process on $T_m$ with parameter $1/\ln
m$.
\end{lemma}

\subsection{Proof of Theorem \ref{T2}}
\label{Sdestruct-mainproof} 
Lemma \ref{L4} suggests that one first proves finite-dimensional
convergence of generation $1$ in the tree of components and then uses the
branching property to transfer the convergence to higher generations. For
the sake of clarity, we restate the result for the first generation.
\begin{proposition}
\label{P1}
For every $\ell\in\mathbb{N}$, there is the convergence in
distribution as $n\rightarrow\infty$, 
\begin{equation*}
  \left(\left(\frac{\ln
        n}{n}\pB_{\sigma(1)}^{(n)},b_{\sigma(1)}^{(n)}\right),\dots,\left(\frac{\ln
        n}{n}\pB_{\sigma(\ell)}^{(n)},b_{\sigma(\ell)}^{(n)}\right)\right)\Longrightarrow
  \left((\pZ_1,z_1),\dots,(\pZ_\ell,z_\ell)\right). 
\end{equation*}
\end{proposition}
We recall that the sequence $(\pB_i^{(n)}: i\in\mathbb{N})$ can be
identified with the sizes of the tree components which are cut from the 
component containing the root $0$ one after the other in the algorithm for isolating the
root. The coupling of Iksanov and M\"ohle provides information on the
component sizes, but it is no longer valid beyond
the last passage time $L(n)$ of the random walk. Our way out is to show
that for fixed $\ell\in\N$, the probability that only those tree components
which are born before time $t$ do contribute to the $\ell$ largest is as
close to $1$ as we wish for $t$ and $n$ sufficiently large.  For fixed
$t>0$, Lemma 4 of Bertoin \cite{Be1} shows that in the destruction process
up to time $t$, about $(1-\e^{-t})n/\ln n$ edges have been removed from the
root component. The first statement of \eqref{EQLn} then implies that we
are in a regime where the coupling applies.

Corollary $2$ of \cite{Be1} already provides us with a limit result for
the tree components that are cut off from
the root component before time $t(n)>0$. Let
$$
\left(\pB_{i,t(n)},b_{i,t(n)}\right) = \left\{\begin{array}{l@{\quad \mbox{for\ }}l}
    (\pB_i^{(n)},b_i^{(n)}) & b_i^{(n)}<t(n)\\
    (0,\infty) & b_i^{(n)}\geq t(n)\end{array}\right.,
$$
and write $$\pC_{\emptyset,t(n)} = n+1-\sum_{j=1}^\infty\pB_{j,t(n)}$$ for
the size of the root component of $T_n$ at time $t(n)$. The parameter $n$
is dropped here as a superscript, since it already appears in $t(n)$.
Furthermore, denote by $\mu=\mu(n,t(n))$ a random permutation that sort the
elements $(\pB_{i,t(n)},b_{i,t(n)})$ in the decreasing order of their first
coordinate, i.e. such that
$$\pB_{\mu(1),t(n)}\geq \pB_{\mu(2),t(n)}\geq\dots.$$

\begin{lemma}{\rm (Bertoin \cite{Be1})}
\label{L5}
Let $(t(n) : n\in\N)$ be a sequence of times converging to some
$t>0$. Then, for each $\ell\in\mathbb{N}$, there is the convergence in
distribution as $n\rightarrow\infty$,
\begin{equation*}
  \left(\left(\frac{\ln
        n}{n}\pB_{\mu(1),t(n)},b_{\mu(1),t(n)}\right),\dots,\left(\frac{\ln
        n}{n}\pB_{\mu(\ell),t(n)},b_{\mu(\ell),t(n)}\right)\right)\Longrightarrow
  \left((\pZ_{1,t},z_{1,t}),\dots,(\pZ_{\ell,t},z_{\ell,t})\right),
\end{equation*}
where $(\pZ_{1,t},z_{1,t}), (\pZ_{2,t},z_{2,t}),\dots$ are the atoms of a Poisson
random measure on $(0,\infty)\times (0,t)$ with intensity $a^{-2}\dt a
\otimes \e^{-s}\dt s$, ranked in the decreasing order of their first
coordinate.  Moreover, $\pC_{\emptyset,t(n)}\sim e^{-t}n$ in
probability as $n\rightarrow\infty$.
\end{lemma}
\noindent{\bf Remark.}
Basic properties of Poisson random measures show that the sequence of atoms
$((\pZ_{i,t},z_{i,t}) : i\in\N)$ can be obtained from $((\pZ_i,z_i):i\in \N)$
by ``squeezing-out'' the elements $(\pZ_i,z_i)$ with $z_i\geq t$.
Formally, conditionally on $(z_i:i\in \N)$, define a map $\gamma :
\N_0\rightarrow \N_0$ by setting $\gamma(0) = 0$, and then
for $i=1,2,\dots$, $\gamma(i)=\inf\{j> \gamma(i-1): z_{j}< t\}$. Then the
sequence $(\pZ_{\gamma(i)},z_{\gamma(i)}) : i\in\N)$ has the same
distribution as $((\pZ_{i,t},z_{i,t}) : i\in\N)$. This point of view is useful
for the proof of Proposition \ref{P1}, which we give now.

\proofof{\bf Proposition \ref{P1}:}\hskip10pt We will reduce the statement
to Lemma \ref{L5} by showing that the $\ell$ largest tree components of
generation $1$ are with high probability produced before time $t$, provided
$t$ and $n$ are sufficiently large. We fix $\ell\in\mathbb{N}$, $\ve>0$ and let
$f:([0,\infty)\times[0,\infty))^\ell\rightarrow[0,1]$ be a continuous
function. Recall that $(z_j)_{j\in\mathbb{N}}$ is a family of
i.i.d. standard exponentials, which is independent of the family
$(\pZ_j)_{j\in\mathbb{N}}$. Choosing $t$ so large such that
$\P\left(\max\{z_1,\dots,z_\ell\} >t\right)\leq \ve$, we obtain by the
remark above
$$
\left|\E\left[f((\pZ_1,z_1),\dots,(\pZ_\ell,z_\ell))\right]-\E\left[f((\pZ_{1,t},z_{1,t}),\dots,(\pZ_{\ell,t},z_{\ell,t}))\right]\right|
\leq \ve.
$$
We will now prove that if $t$ is large, then for all $n$ large enough also
\begin{align}\label{C2E1}
  \lefteqn{\left|\E\left[f\left(\left(\frac{\ln
              n}{n}\pB_{\sigma(1)}^{(n)},b_{\sigma(1)}^{(n)}\right),\dots,\left(\frac{\ln
              n}{n}\pB_{\sigma(\ell)}^{(n)},b_{\sigma(\ell)}^{(n)}\right)\right)\right]\right.}\nonumber\\
  &-\left.\E\left[f\left(\left(\frac{\ln
            n}{n}\pB_{\mu(1),t}^{(n)},b_{\mu(1),t}^{(n)}\right),\dots,\left(\frac{\ln
            n}{n}\pB_{\mu(\ell),t}^{(n)},b_{\mu(\ell),t}^{(n)}\right)\right)\right]\right| 
  \leq \ve.
\end{align}
Here, since $t$ is fixed, we write $\pB_{\mu(\ell),t}^{(n)}$ instead of $\pB_{\mu(\ell),t(n)}$
and similarly for the birth times. First, it follows from Lemma \ref{L5} that for $t_0>0$, there exists $\delta>0$ such
that for each $t>t_0$ and for all $n$ sufficiently large,
$$\P\left(\pB_{\mu(\ell),t}^{(n)}\geq\delta n/\ln n\right)\geq 1-\ve/2.$$ Next, if
$\pB_{\ast,\geq t}^{(n)}$ is the size of the largest tree component amongst those which were cut off from the root component in the
destruction process on $T_n$ at a time $\geq t$,
then on the event $\{\pB_{\mu(\ell),t}^{(n)}\geq\delta n/\ln n\}\cap
\{\pB_{\ast,\geq t}^{(n)}<\delta n/\ln n\}$, there is the equality of random vectors
$$
\left(\left(\pB_{\mu(1),t}^{(n)},b_{\mu(1),t}^{(n)}\right),\dots,\left(\pB_{\mu(\ell),t}^{(n)},b_{\mu(\ell),t}^{(n)}\right)\right)
=
\left(\left(\pB_{\sigma(1)}^{(n)},b_{\sigma(1)}^{(n)}\right),\dots,\left(\pB_{\sigma(\ell)}^{(n)},b_{\sigma(\ell)}^{(n)}\right)\right).
$$
Therefore, \eqref{C2E1} follows if we show that for large $t$ and all large
$n$,
\begin{equation}
\label{C2E2}
\P\left(\pB_{\ast,\geq t}^{(n)}<\delta n/\ln n\right)\geq 1-\ve/2.
\end{equation}
Write $m=m(t,n)$ for the number of edges of the root component in the
destruction process at time $t$. By the splitting property, conditionally
on $m$, the variable $\pB_{\ast,\geq t}^{(n)}$ is distributed as the size of the
largest tree component which was produced by the algorithm for isolating
the root of a RRT of size $m+1$. We now claim that $(\ln m/m)\pB_{\ast,\geq t}^{(n)}$
converges in distribution as $m\rightarrow\infty$ to the largest atom of
a Poisson random measure on $(0,\infty)$ with intensity $a^{-2}\dt a$.

Indeed, if $\xi_1,\xi_2,\dots$ is a sequence of of i.i.d. copies of $\xi$,
see \eqref{Estepdistr}, then for $a>0$, the number of indices $j\leq k$
such that $\xi_j>am/\ln m$ is binomially distributed with parameters $k$
and $\lceil am/\ln m\rceil^{-1}$. Combining the first part of~\eqref{EQLn}
with Theorem 16.16 of Kallenberg \cite{Ka}, we deduce that the largest
variable among $\xi_1,\dots,\xi_{L(m)}$, normalized by a factor $\ln m/m$,
converges in distribution to largest atom of a Poisson random measure on
$(0,\infty)$ with intensity $\nu(da)=a^{-2}\dt a$. Clearly, under the
coupling of Iksanov and M\"ohle, $m+1-S_{L(m)}$ is the size of the
remaining root component after $L(m)$ edge removals in the algorithm for
isolating the root. Since $m+1-S_{L(m)} = o(m/\ln m)$ in probability by
\eqref{EQLn}, an appeal to the coupling proves our claim about
$\pB_{\ast,\geq t}^{(n)}$.

We finally notice that by the second part of Lemma \ref{L5}, $m+1 =
\pC_{\emptyset,t}^{(n)}\sim \e^{-t}n$ in probability, that is $(\ln
n/n)\pB_{\ast,\geq t}^{(n)}$ converges in distribution to the largest atom of a
Poisson random measure on $(0,\infty)$ with intensity $\e^{-t}a^{-2}\dt
a$. Choosing $t=t(\delta)$ large enough, \eqref{C2E2} follows. Since
$\ve>0$ can be chosen arbitrarily small, an application of the triangle
inequality together with Lemma \ref{L5} finishes the proof of Proposition
\ref{P1}.  \QED

Now we are in position to prove Theorem \ref{T2}. The line of argumentation
is similar to that in the proof of Proposition \ref{Ptreeofsubtreesizes}.

\proofof{\bf Theorem \ref{T2}:}\hskip10pt The convergence of
$((1/n)\pB_{\emptyset}^{(n)},b_{\emptyset}^{(n)})$ is
trivial, and Proposition \ref{P1} shows the convergence of generation
$1$. Let us now show that also
\begin{equation}
\label{eq:2dlaws}
\left(\left(\frac{(\ln n)^{|u|}}{n}\pB_{\sigma(u)}^{(n)},b_{\sigma(u)}^{(n)}\right):u\in\pU,|u|\leq 2\right)\Longrightarrow ((\pZ_u, z_u):u\in\pU,|u|\leq 2)
\end{equation} 
as $n\rightarrow\infty$ in the sense of finite-dimensional laws. Let
$\ell\in\mathbb{N}$. Employing Lemma \ref{L4} and Proposition \ref{P1}, it
suffices to show that for
$g,f^{(i)}:[0,\infty)\times[0,\infty]\rightarrow [0,1]$ bounded and
uniformly continuous, and $j,\ell\in\N$,
\begin{align*}
  \lefteqn{\E\left[g\left(\frac{\ln
          n}{n}\pB_{\sigma(j)}^{(n)},b_{\sigma(j)}^{(n)}\right)f^{(1)}\left(\frac{\ln^2
          n}{n}
\pB_{\sigma(j1)}^{(n)},b_{\sigma(j1)}^{(n)}\right)\dots
      f^{(\ell)}\left(\frac{\ln^2 n}{n}\pB_{\sigma(j\ell)}^{(n)},b_{\sigma(j\ell)}^{(n)}\right)\right]}\\
  &\rightarrow  \E\left[g(\pZ_j,z_j)f^{(1)}(\pZ_ja_{1},z_j +b_{1})\dots 
    f^{(\ell)}(\pZ_ja_{\ell},z_j+b_{\ell})\right],\hspace{3cm}
\end{align*} 
where for $i=1,\dots,\ell$, 
$(a_{i},b_{i})$ is the atom with the $i$th largest first coordinate of a
Poisson random measure on $(0,\infty)\times(0,\infty)$ with intensity
$a^{-2}\dt a\otimes \e^{-r}\dt r$. We consider only the case $\ell=1$.  By
Lemma \ref{L5}, we have for each integer $m$ with
$\P(\pB_{\sigma(j)}^{(n)}=m)>0$ and almost all $s>0$ the equality of the
conditional densities
\begin{equation*}
 \E\left[f^{(1)}\left(\frac{\ln^2n}{n}\pB_{\sigma(j1)}^{(n)},
        b_{\sigma(j1)}^{(n)}\right)\Big|\,\left(\pB_{\sigma(j)}^{(n)},b_{\sigma(j)}^{(n)}\right)=(m,s)\right] =
  \E_m\left[f^{(1)}\left(\frac{\ln^2
        n}{n}\pB_{\ast}^{(m)},s+\frac{\ln n}{\ln m}b_{\ast}^{(m)}\right)\right], 
\end{equation*}
where $\E_m$ is the mathematical expectation starting from a random
recursive tree with $m$ vertices, and under $\E_m$, 
$(\pB_{\ast}^{(m)},b_{\ast}^{(m)})$ is in the first coordinate
the size and in the second the birth time of the largest tree component of the
first generation produced by a destruction process on $T_m$ with parameter
$1/\ln m$. Now if $m\sim (n/\ln n)a$ for some fixed $a>0$, $m=m(n)$
integer-valued, we obtain from Proposition \ref{P1} that
\begin{equation*}
  \E_m\left[f^{(1)}\left(\frac{\ln^2n}{n}\pB_{\ast}^{(m)},s+\frac{\ln
        n}{\ln m}b_{\ast}^{(m)}\right)\right]\sim \E\left[f^{(1)}(aa_1,s+b_1)\right],
\end{equation*} 
where $(a_1,b_1)$ is the atom with the largest first coordinate of a
Poisson random measure on $(0,\infty)\times(0,\infty)$ with intensity
$a^{-2}\dt a\otimes \e^{-r}\dt r$. On the other hand, we already know that
the pair $(\frac{\ln n}{n}\pB_{\sigma(j)}^{(n)},b_{\sigma(j)}^{(n)})$
converges in distribution as $n\rightarrow \infty$ towards
$(\pZ_j,z_j)$. Since the map
\begin{equation*}
  ((a,b),(a_1,b_1))\mapsto g(a,b)f^{(1)}(aa_1,b+b_1)
\end{equation*}
is uniformly continuous on bounded sets, this
establishes~\eqref{eq:2dlaws}. The arguments can now easily be extended to
the subsequent generations, and the theorem is proved.\QED

\subsection{Applications of the destruction process and remarks}
\label{Sapps}

\subsubsection{Connection to Bernoulli bond percolation}
In \cite{Be1}, Bertoin uses the coupling of Iksanov and M\"ohle to study
the asymptotic sizes of the largest and next largest percolation clusters
of a supercritical Bernoulli bond percolation on $T_n$ with
parameter 
\begin{equation}
\label{pnsuperc}
p(n)= 1-t/\ln n +o(1/\ln n),\quad t>0\hbox{ fixed.}
\end{equation} 
Let us recall his strategy. If the destruction process (with parameter
$1/\ln n$) is stopped at time $t(n) = -\ln n\times\ln p(n)$,
then one observes a Bernoulli bond percolation on $T_n$ with parameter
$p(n)$.  Under this coupling, the tree components born in the destruction
process up to time $t(n)$ contain the non-root percolation clusters of $T_n$. In
fact, each such percolation cluster of $T_n$ can be identified with a subtree of
a tree component rooted at the same vertex, meaning that within its
surrounding component, the percolation cluster forms the root cluster.

The usefulness of this point of view comes from two facts. Firstly, we know
from the second part of Lemma \ref{L5} that in the regime \eqref{pnsuperc}, the
root cluster of a RRT $T_m$ has size $\sim \e^{-t}m$ as
$m\rightarrow\infty$. Secondly, the asymptotic sizes of the tree components
can be specified (see Proposition \ref{P2}). In order to reveal the inner
root percolation cluster inside a tree component, the latter has to be
``unfrozen'', i.e. some additional edges have to be erased. This approach
was used by Bertoin \cite{Be1} to study the sizes of the root percolation
clusters inside the tree components of the first generation, and our aim is
to outline how these ideas can be extended to all clusters. We first lift
the convergence of Lemma \ref{L5} to higher generations. Towards this end,
let
$$
\left(\pB_{u,t(n)},b_{u,t(n)}\right) = \left\{\begin{array}{l@{\quad \mbox{for\ }}l}
    (\pB_u^{(n)},b_u^{(n)}) & b_u^{(n)}<t(n)\\
    (0,\infty) & b_u^{(n)}\geq t(n)\end{array}\right..
$$
Then we can use Lemma \ref{L5} instead of Proposition \ref{P1} to obtain a
limit result for the ranked version
$((\pB_{\sigma(u),t(n)},b_{\sigma(u),t(n)}) : u\in\pU)$. Here, by a small
abuse of notation, $\sigma:\pU\rightarrow\pU$
is a random bijection that sorts the children of each element
$(\pB_{u,t(n)},b_{u,t(n)})$ in the decreasing order of
their first coordinate, keeping the parent-child
relation. The limit process $((\pZ_{u,t}, z_{u,t}) : u\in\pU)$ is obtained
from $((\pZ_u, z_u) : u\in\pU)$ by first removing those pairs $(\pZ_u,z_u)$
with $z_u\geq t$ and then by a relabeling of the remaining
elements. Alternatively, in accordance with Lemma \ref{L5}, the law of the
limit can also be specified as follows.
\begin{enumerate}
 \item $(\pZ_{\emptyset,t},z_{\emptyset,t}) = (1,0)$ almost surely;
 \item for every $k=0,1,2,\dots,$ conditionally on $((\pZ_{v,t},z_{v,t}):
   v\in\pU,|v|\leq k)$, the sequences $((\pZ_{uj,t},z_{uj,t}))_{j\in\N}$
   for the vertices $u\in\pU$ at generation $|u|=k$ are independent, and
   each sequence $((\pZ_{uj,t},z_{uj,t}-z_{u,t}))_{j\in\N}$ is distributed as
   the family of the atoms of a Poisson random measure on $(0,\infty)\times
   (0,t-z_{u,t})$ with intensity $\pZ_{u,t}a^{-2}\dt a\otimes \e^{-r}\dt r$,
   ranked in the decreasing order of the first coordinate.
 \end{enumerate}
The analog of Theorem \ref{T2} for the tree components born up to time
 $t(n)$ then reads as follows.
\begin{proposition}
\label{P2}
As $n\rightarrow\infty$, in the sense of finite-dimensional distributions,
$$
\left(\left(\frac{(\ln n)^{|u|}}{n}\pB_{\sigma(u),t(n)},
    b_{\sigma(u),t(n)}\right) : u\in\pU\right)\Longrightarrow
  \left(\left(\pZ_{u,t}, z_{u,t}\right) : u\in\pU\right).
$$
\end{proposition}
Now the tree components have to be unfrozen to observe the percolation
clusters inside. Write $\tau_u^{(n)}$ for the tree component whose size and
birth time is stored in $(\pB_{\sigma(u),t(n)},b_{\sigma(u),t(n)})$
(with $\tau_u^{(n)}=\emptyset$ if there is no such component, and $\tau_{\emptyset}^{(n)}=T_n$). Say we want
to determine the size of the root percolation cluster inside the tree component $\tau_u^{(n)}$. This
component was cut off from a bigger subtree at time
$b=b_{\sigma(u),t(n)}$. By the memoryless property of exponential
variables, we are therefore lead to perform a Bernoulli bond percolation on
$\tau_u^{(n)}$ with parameter $\exp(-(t(n)-b)/\ln n)$, and adapting the
arguments of \cite{Be1}, we deduce that the root cluster $c^{(n)}_u$ of
$\tau_u^{(n)}$ has size 
$$|c^{(n)}_u|\sim \e^{-(t-b)}\pB_{\sigma(u),t(n)}.$$
More generally, denote by $c^{(n)}_u$ for $u\in\pU$ the percolation cluster
with the same root as $\tau_u^{(n)}$ (under our coupling with the
destruction process). In the percolation regime $1-p(n)\sim t/\ln n$, we
have
$$(1-p(n))^{-k}n^{-p(n)}\sim t^{-k}\e^t(\ln n)^k n^{-1}.$$ Using Proposition
\ref{P2}, the last two displays and similar arguments as in the proof of Theorem
\ref{T2}, we obtain the following limit result for the cluster sizes $|c^{(n)}_u|$.
\begin{corollary}
  \label{C2}
 As $n\rightarrow\infty$, in the sense of finite-dimensional distributions,
$$
\left(\frac{(1-p(n))^{-|u|}}{n^{p(n)}}|c^{(n)}_u| : u\in\pU\right)\Longrightarrow
\left(t^{-|u|}\exp\left(z_{u,t}\right)\pZ_{u,t}: u\in\pU\right).
$$
\end{corollary}
Extending the arguments of \cite[Lemma 6, 7]{Be1} to higher levels in the
tree, we moreover  see that Corollary \ref{C2} remains
true if we apply our usual ranking operation to both sides. Denote by
$\pCf^{(n)}=(\pCf^{(n)}_u: u\in\pU)$ the ranked version of $(|c^{(n)}_u|:
u\in\pU)$, i.e. $\pCf^{(n)}_u=|c^{(n)}_{\tilde{\sigma}(u)}|$, where
$\tilde{\sigma}:\pU\rightarrow\pU$ is a random bijection sorting the
children $(|c^{(n)}_{ui}|:i\in \N)$ of each element $|c^{(n)}_u|$ in the decreasing order, such that the
parent-child relation is preserved. For the right hand side, let us write
$\pG_u =t^{-|u|}\exp(z_{u,t})\pZ_{u,t}$ and $(\pG_{\theta(u)} : u\in\pU)$
for the ranked version of $(\pG_u:u\in\pU)$. Then the convergence in  
Corollary \ref{C2} transfers to the ranked versions, i.e. 
$$\left(\frac{(1-p(n))^{-|u|}}{n^{p(n)}}\pCf^{(n)}_u:
  u\in\pU\right)\Longrightarrow (\pG_{\theta(u)}  : u\in\pU)$$ 
in the sense of finite-dimensional distributions in the regime \eqref{pnsuperc}. It is
now instructive to compare this last convergence result with Theorem
\ref{T1}.

We first remark that as for the tree of cluster sizes $\pC^{(n)}$ from Section
\ref{Smainresults}, the process $(\pCf^{(n)}_u: u\in\pU)$ stores the size
of every percolation cluster of $T_n$. Both $\pC^{(n)}_{\emptyset}$ and
$\pCf^{(n)}_{\emptyset}$ encode the size of the cluster containing
$0$. But besides that, the two encodings are different.
Most importantly, if we look at some specific percolation cluster of $T_n$
and ask for the vertex $u\in\pU$ to which the size of this cluster is
attached in the process $(\pCf^{(n)}_u : u\in\pU)$, we observe that its
level $|u|=k$ does not merely depend on the total number of removed edges
which separate the cluster from the vertex $0$, but also on the order in which these
edges were removed.

To stress the difference in the encodings, call a percolation
cluster encoded by some $c^{(n)}_u$ with $|u|=k$ a cluster of rank $k$. In
terms of our classification of clusters into generations from Section
\ref{Smainresults}, a cluster of generation $k\geq 1$ with root node $v$ can be a cluster of
rank $1\leq \ell\leq k$; the rank depends on the order in which the $k$
erased edges on the path from $0$ to $v$ were removed in the
destruction process. Conversely, a cluster of rank $\ell$ with root node
$v$ can be a cluster of generation $k$ for $\ell\leq k\leq\, $dist$(0,v)$,
where dist$(\cdot,\cdot)$ denotes the graph distance on $T_n$ before the
percolation was performed.
 
Figure $3$ illustrates the difference in the encoding by $\pC^{(n)}$ and
$\pCf^{(n)}$, respectively. We tacitly assume that the tree of cluster
sizes $\pC^{(n)}$ is defined in terms of the final state of a percolation
on $T_n$ which is used to define $\pCf^{(n)}$. For example, the
cluster $\{5,9\}$ is a cluster of rank $1$, since the edge joining $5$ to
its parent $3$ was the first edge from the path connecting $0$ to $5$ which was
removed in the destruction process. On the other hand, $\{5,9\}$ is a
cluster of generation $2$, since it is disconnected from $0$ by two
deleted edges in the final outcome of percolation.

\begin{figure}[ht]
\begin{center}\parbox{5cm}{\includegraphics[width=4.5cm]{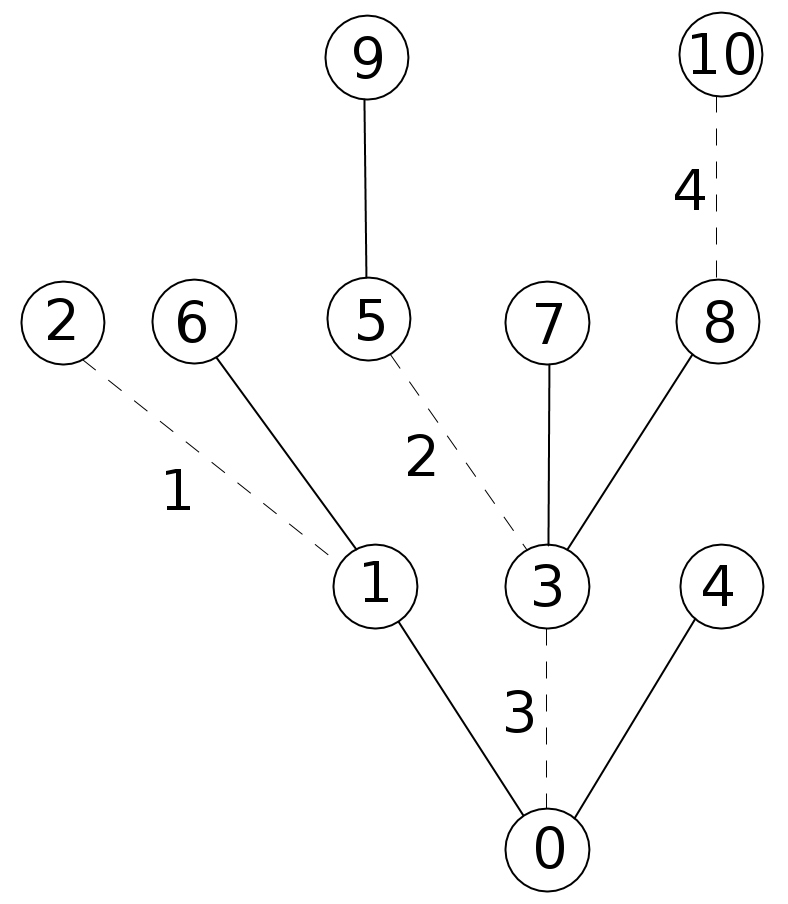}}
\parbox{5.5cm}{\includegraphics[width=5.5cm]{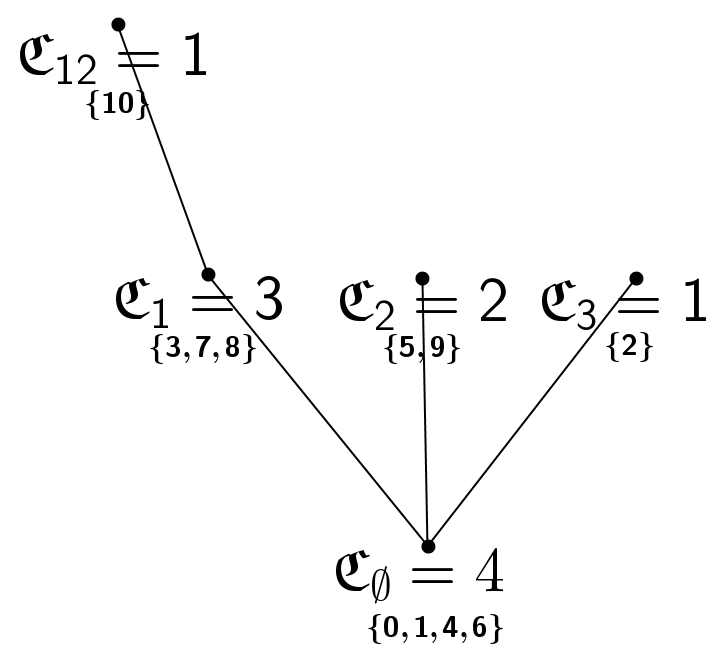}}
\parbox{6.0cm}{\includegraphics[width=6.0cm]{rrt-treeofclusters-nolabels.png}}
\end{center}
\centerline{\bf Figure 3} \centering{\sl Left: Percolation on a
  RRT with vertex labels $0,1,\dots,10$. The erased edges are indicated
  by dashed lines, and their labels indicate the order in which the edges were removed.}\\
{\sl Middle: The cluster encoding by $\pCf^{(n)}$. Note that several orderings of
  edge removals give rise to the same tree.}\\
{\sl Right: The tree of cluster sizes $\pC^{(n)}$ defined in Section \ref{Smainresults}.}
\end{figure}
Recall the description of $((\pZ_{u,t},z_{u,t}) :u\in\pU)$ from above.
We now observe that conditionally on $(\pZ_{u,t},z_{u,t})$, the family
$(\pG_{\theta(uj)} : j\in\mathbb{N})$ is distributed as the sequence $b_1>
b_2> \dots$ of the atoms of a Poisson random measure on $(0,\infty)$ with
intensity
$$t^{-(|u|+1)}(t-z_{u,t})\exp(z_{u,t})\pZ_{u,t}a^{-2}\dt a.$$ Indeed,
$\pG_{uj}=t^{-(|u|+1)}\exp(z_{u,t})\exp(z_{uj,t}-z_{u,t})\pZ_{uj,t}$, and
given $(\pZ_{u,t}, z_{u,t})$, the image of the measure $\pZ_{u,t}a^{-2}\dt
a\otimes \e^{-r}\dt r$ on $(0,\infty)\times (0,t-z_{u,t})$ by the map
$(a,s)\mapsto t^{-(|u|+1)}\exp(z_{u,t})\exp(s)a$ is
$t^{-(|u|+1)}(t-z_{u,t})\exp(z_{u,t})\pZ_{u,t}a^{-2}\dt a$ on $(0,\infty)$.

Since $(Z_{\emptyset,t},z_{\emptyset,t}) =(1,0)$, we deduce from this
characterization that the sequences $(\pG_{\theta(j)}:j\in\N)$ and $(\pZ_j
:j\in\N)$ have the same distribution, which implies that the
finite-dimensional limits of $(\pC_u^{(n)}: u\in\pU, |u|\leq 1)$ and
$(\pCf_u^{(n)} : u\in\pU, |u|\leq 1)$ agree (under our normalizations). 

In fact, this already follows from our previous considerations: We have
seen in the proof of Corollary \ref{C1} that the largest non-root clusters
are of generation $1$, and every such cluster is necessarily a cluster of
rank $1$ (but not every cluster of rank $1$ is of generation $1$, see
cluster $\{5,9\}$ in Figure $3$). For higher levels in the trees
$\pC^{(n)}$ and $\pCf^{(n)}$, the limits do however not agree. This comes
from the fact that clusters of generation $k\geq 2$ can represent clusters
of a strictly lower rank $1\leq \ell < k$. Roughly speaking, if such a
cluster has a size of order $(1-p(n))^{k}n^{p(n)}$, it is visible in the
limit under the encoding by $\pC^{(n)}$, while it is not under the encoding
by $\pCf^{(n)}$.

\subsubsection{Connection to the cut-tree}
The tree of components is related to the so-called {\it cut-tree}, which is
defined in terms of a discrete-time destruction process, where edges are
removed according to some order, for example a random uniform order. 

\begin{figure}[ht]
\begin{center}\parbox{8cm}{\includegraphics[width=6.5cm]{rrt-destruction.png}}
  \parbox{6.5cm}{\includegraphics[width=6.5cm]{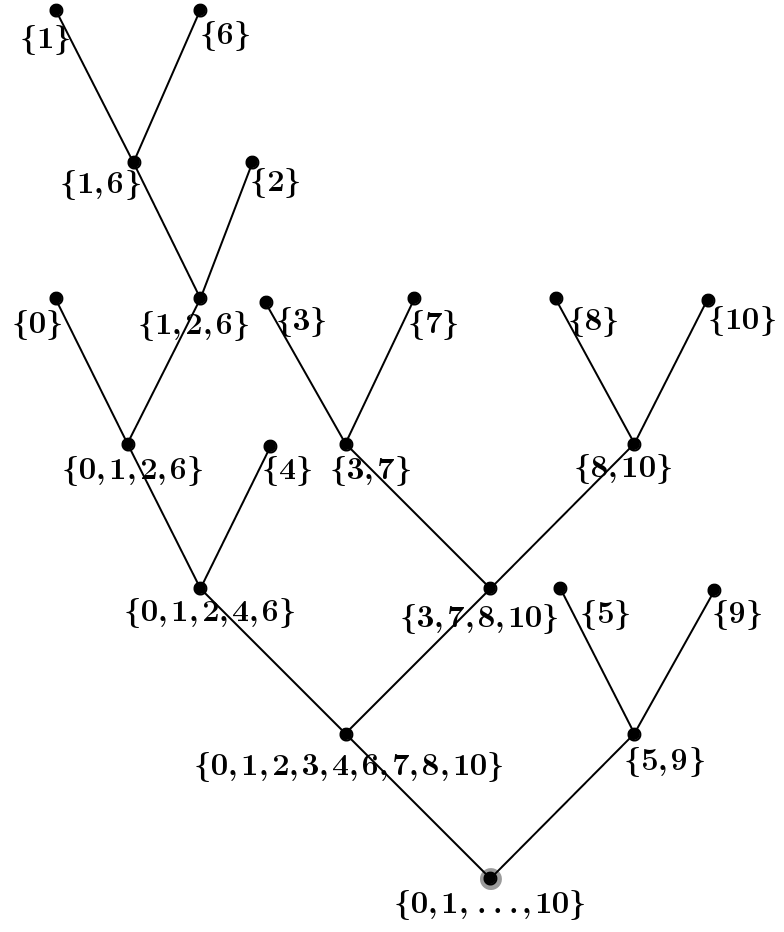}}
\end{center}
\centerline{\bf Figure 4} \centering{\sl Left: The same recursive tree as in Figure 2, with the same order of edge removals.}\\
{\sl Right: The corresponding cut-tree.}
\end{figure}

More specifically, the cut-tree is a rooted binary tree which encodes the
destruction of a tree $T$ on a finite vertex set $V$ in the following
way. The root vertex is given by the set $V$. Then, if the first edge is
removed, $T$ splits into two subtrees with respective vertex sets
$V_1$ and $V_2$, and these vertex sets are attached as the two children to
the root $V$. The construction is then iterated in the natural way - if,
for example, the next edge is removed from the subtree with vertex
set $V_1$, the latter splits into two vertex sets $V_{1,1}$ and $V_{1,2}$,
which are regarded as the two children of $V_1$. In particular, the leaves
of the cut-tree can be identified with the vertices of $T$.  

Unlike the tree of components, the cut-tree stores the vertex sets of the
tree components and not merely their sizes. For example, in Figure 4 the
vertex sets of the tree components of the first generation, i.e. $\{5,9\}$,
$\{3,7,8,10\}$, $\{4\}$ and $\{1,2,6\}$ (in the order of their appearance),
are represented by the vertices which are attached to the branch from the
root $\{0,1,\dots,10\}$ to the leaf $\{0\}$.

The cut-tree has been analyzed for Cayley trees and random recursive trees
by Bertoin in \cite{Be0} and \cite{Be3}, and then by Bertoin and Miermont
\cite{BM} and Dieuleveut \cite{Di} for Galton-Watson trees.  Their results
can be used to obtain limit theorems for the number of steps to isolate a
certain family of nodes, and in a similar direction, we believe that the
tree of components can prove helpful, too.

\vspace{1cm}
\noindent{\bf Acknowledgments.} I would like to thank Jean Bertoin for his
generosity in sharing his insights, and for helpful discussions.

\end{document}